\documentclass[reqno]{amsart}

\usepackage{tikz-cd}

\usepackage{fullpage,dsfont}

\usepackage{hyperref} 

\usepackage{parskip}
\makeatletter 
\def\thm@space@setup{%
 \thm@preskip=\parskip \thm@postskip=0pt
}
\def\th@remark{%
  \thm@headfont{\itshape}%
  \normalfont 
  \thm@preskip\parskip \thm@postskip=0pt
}
\makeatother

\usepackage[nobysame,alphabetic,initials,msc-links]{amsrefs}

\usepackage[final]{pdfpages}

\usepackage{multirow}

\usepackage{array}

\usepackage{tikz-cd}
\usetikzlibrary{knots}
\usepackage{spath3}

\usepackage{float}

\DefineSimpleKey{bib}{how}
\DefineSimpleKey{bib}{mrclass}
\DefineSimpleKey{bib}{mrnumber}
\DefineSimpleKey{bib}{fjournal}
\DefineSimpleKey{bib}{mrreviewer}

\renewcommand{\PrintDOI}[1]{%
  \href{http://dx.doi.org/#1}{{\tt DOI:#1}}%
}
\renewcommand{\eprint}[1]{#1}
\BibSpec{book}{%
    +{}  {\PrintPrimary}                {transition}
    +{.} { \PrintDate}                  {date}
    +{.} { \textit}                     {title}
    +{.} { }                            {part}
    +{:} { \textit}                     {subtitle}
    +{,} { \PrintEdition}               {edition}
    +{}  { \PrintEditorsB}              {editor}
    +{,} { \PrintTranslatorsC}          {translator}
    +{,} { \PrintContributions}         {contribution}
    +{,} { }                            {series}
    +{,} { \voltext}                    {volume}
    +{,} { }                            {publisher}
    +{,} { }                            {organization}
    +{,} { }                            {address}
    +{,} { }                            {status}
    +{,} { \PrintDOI}                   {doi}
    +{,} { \PrintISBNs}                 {isbn}
    +{}  { \parenthesize}               {language}
    +{}  { \PrintTranslation}           {translation}
    +{;} { \PrintReprint}               {reprint}
    +{.} { }                            {note}
    +{.} {}                             {transition}
    +{}  {\SentenceSpace \PrintReviews} {review}
}
\BibSpec{article}{%
    +{}  {\PrintAuthors}                {author}
    +{,} { \textit}                     {title}
    +{.} { }                            {part}
    +{:} { \textit}                     {subtitle}
    +{,} { \PrintContributions}         {contribution}
    +{.} { \PrintPartials}              {partial}
    +{,} { }                            {journal}
    +{}  { \textbf}                     {volume}
    +{}  { \PrintDatePV}                {date}
    +{,} { \issuetext}                  {number}
    +{,} { \eprintpages}                {pages}
    +{,} { }                            {status}
    +{,} { \PrintDOI}                   {doi}
    +{,} { \eprint}        {eprint}
    +{}  { \parenthesize}               {language}
    +{}  { \PrintTranslation}           {translation}
    +{;} { \PrintReprint}               {reprint}
    +{.} { }                            {note}
    +{.} {}                             {transition}
    +{}  {\SentenceSpace \PrintReviews} {review}
}
\BibSpec{collection.article}{%
    +{}  {\PrintAuthors}                {author}
    +{,} { \textit}                     {title}
    +{.} { }                            {part}
    +{:} { \textit}                     {subtitle}
    +{,} { \PrintContributions}         {contribution}
    +{,} { \PrintConference}            {conference}
    +{}  {\PrintBook}                   {book}
    +{,} { }                            {booktitle}
    +{,} { \PrintDateB}                 {date}
    +{,} { pp.~}                        {pages}
    +{,} { }                            {publisher}
    +{,} { }                            {organization}
    +{,} { }                            {address}
    +{,} { }                            {status}
    +{,} { \PrintDOI}                   {doi}
    +{,} { \eprint}        {eprint}
    +{}  { \parenthesize}               {language}
    +{}  { \PrintTranslation}           {translation}
    +{;} { \PrintReprint}               {reprint}
    +{.} { }                            {note}
    +{.} {}                             {transition}
    +{}  {\SentenceSpace \PrintReviews} {review}
}
\BibSpec{misc}{%
  +{}{\PrintAuthors}  {author}
  +{,}{ \textit}      {title}
  +{.}{ }             {how}
  +{}{ \parenthesize} {date}
  +{,} { available at \eprint}        {eprint}
  +{,}{ available at \url}{url}
  +{,}{ }             {note}
  +{.}{}              {transition}
}

\usepackage{amssymb, amsfonts, amsxtra, amsmath}
\usepackage{mathrsfs}
\usepackage{mathdots}
\usepackage{wasysym}
\usepackage[all]{xy}
\usepackage{bbm}
\usepackage{calc}
\usepackage{accents}

\numberwithin{equation}{section}

\newtheorem{Theorem}{Theorem}[section]
\newtheorem*{Theorem*}{Theorem}
\newtheorem{Def}[Theorem]{Definition}
\newtheorem{Lem}[Theorem]{Lemma}
\newtheorem{Prop}[Theorem]{Proposition}

\newtheorem{Rem}[Theorem]{Remark}

\newtheorem{Exa}[Theorem]{Example}

\newcommand\bp{\begin{proof}}
\newcommand\ep{\end{proof}}

\mathchardef\mhyph="2D

\DeclareMathOperator{\fin}{\mathrm{fin}}
\DeclareMathOperator{\loc}{\mathrm{loc}}

\DeclareMathOperator{\Hilb}{\mathrm{Hilb}}

\DeclareMathOperator{\id}{\mathrm{id}}

\DeclareMathOperator{\Rep}{\mathrm{Rep}}
\DeclareMathOperator{\rd}{\mathrm{d}\!}
\DeclareMathOperator{\Tr}{\mathrm{Tr}}

\DeclareMathOperator{\adm}{\mathrm{adm}}

\DeclareMathOperator{\Ind}{\mathrm{Ind}}
\DeclareMathOperator{\Irr}{\mathrm{Irr}}

\DeclareMathOperator{\Ker}{\mathrm{Ker}}

\DeclareMathOperator{\Lin}{\mathrm{Lin}}

\DeclareMathOperator{\Mod}{\mathrm{Mod}}

\DeclareMathOperator{\Heis}{\mathrm{Heis}}

\newcommand{\ev}{\mathrm{ev}}

\newcommand{\reg}{\mathrm{reg}}

\newcommand{\msG}{\mathscr{G}}

\newcommand{\msI}{\mathscr{I}}

\newcommand{\msMod}{\mathrm{Mod}}
\newcommand{\msRep}{\mathrm{Rep}}

\newcommand{\msP}{\mathscr{P}}

\newcommand{\msU}{\mathscr{U}}

\newcommand{\mfa}{\mathfrak{a}}

\newcommand{\mfg}{\mathfrak{g}}

\newcommand{\mfk}{\mathfrak{k}}
\newcommand{\mfl}{\mathfrak{l}}

\newcommand{\mfn}{\mathfrak{n}}

\newcommand{\mfsl}{\mathfrak{sl}}
\newcommand{\mfso}{\mathfrak{so}}

\newcommand{\mfsu}{\mathfrak{su}}

\newcommand{\mfu}{\mathfrak{u}}

\newcommand{\mcB}{\mathcal{B}}

\newcommand{\mcD}{\mathcal{D}}
\newcommand{\mcE}{\mathcal{E}}

\newcommand{\mcH}{\mathcal{H}}
\newcommand{\Hsp}{\mathcal{H}}

\newcommand{\mcI}{\mathcal{I}}

\newcommand{\mcL}{\mathcal{L}}

\newcommand{\mcO}{\mathcal{O}}

\newcommand{\mcT}{\mathcal{T}}
\newcommand{\mcU}{\mathcal{U}}

\newcommand{\mcX}{\mathcal{X}}

\newcommand{\C}{\mathbb{C}}

\newcommand{\K}{\mathbb{K}}
\newcommand{\N}{\mathbb{N}}
\newcommand{\Q}{\mathbb{Q}}

\newcommand{\R}{\mathbb{R}}

\newcommand{\U}{\mathbb{U}}

\newcommand{\Z}{\mathbb{Z}}

\newcommand{\comp}{\mathrm{comp}}

\newcommand{\opp}{\mathrm{op}}

\newcommand\rhdb{\blacktriangleright}

\newcommand{\ldb}{\{\!\!\{}
\newcommand{\rdb}{\}\!\!\}}
\newcommand{\qdif}[1]{\ldb #1\rdb}

\newcommand{\ldd}{[\![}
\newcommand{\rdd}{]\!]}
\newcommand{\qbra}[1]{\ldd #1\rdd}

\newcommand{\Func}{\mathrm{Func}}

\def\vp{\varphi}
\newcommand{\rphis}[5]{\,_{#1}\vp_{#2} \left( \genfrac{.}{.}{0pt}{}{#3}{#4}
\ ;#5 \right)}

\title{Induction for representations of coideal doubles, with an application to quantum $SL(2,\R)$}
\author{K. De Commer}
\address{Vrije Universiteit Brussel}
\email{kenny.de.commer@vub.be}

\begin{document}
\maketitle

\begin{center}
\emph{In celebration of the 80th birthday of T.H. Koornwinder, with admiration}
\end{center}

\begin{abstract}
We investigate the theory of induction in the setting of doubles of coideal $*$-subalgebras of compact quantum group Hopf $*$-algebras. We then exemplify parts of  this theory in the particular case of quantum $SL(2,\R)$, and compute the decomposition of the regular representation for quantum $SL(2,\R)$ into irreducibles. 
\end{abstract}

\section*{Introduction}

\emph{CQG Hopf $*$-algebras} (or \emph{CQG algebras} for short) were introduced by Koornwinder and Dijkhuizen in \cite{DK94}. They offer a bridge between the purely algebraic theory of Hopf algebras and the analytic theory of compact quantum groups (=CQGs) \cite{Wor87b}. A CQG Hopf $*$-algebra consists of a complex Hopf algebra $(A,\Delta_A)$ with a compatible $*$-structure $*: A \rightarrow A$, that is an anti-linear involution satisfying 
\[
(ab)^* = b^*a^*,\qquad \Delta_A(a^*)= \Delta_A(a)^{*\otimes *},\qquad \forall a,b\in A,
\]
and for which there exists a (necessarily unique) \emph{invariant state} $\Phi_A: A \rightarrow A$, that is, 
\[
\textrm{(state property)}\qquad \Phi_A(1_A) =1,\qquad \Phi_A(a^*a)\geq 0,\qquad \forall a\in A,
\]
\[
\textrm{(invariance property)}\qquad (\Phi_A\otimes \id)\Delta_A(a) = \Phi_A(a)1_A =(\id\otimes \Phi_A)\Delta_A(a),\qquad \forall a\in A. 
\]
For example, any compact (Hausdorff) group $U$ leads to a CQG Hopf $*$-algebra with $A = \mcO(U)$, the $*$-algebra of matrix coefficients of finite-dimensional continuous unitary representations of $U$, with $\Delta$ dual to the product in the group and with, for $\mu$ the (normalized) Haar measure on $U$,
\[
\Phi_A(f) = \int_U f(x) \rd \mu(x),\qquad f\in \mcO(U). 
\] 
In fact, any CQG Hopf $*$-algebra with $A$ commutative arises in this way, for a uniquely determined $U$ (up to isomorphism). For this reason, a general CQG Hopf $*$-algebra is sometimes denoted $\U = (A,\Delta_A)$, and one then writes $A = \mcO(\U)$, referring to $\U$ as the `compact quantum group'. 

From the above picture, it is easy to come up with a notion of a compact quantum subgroup $\K \subseteq \U$: it should consist of a CQG Hopf $*$-algebra $\K = (C,\Delta_C)$ with  a surjective $*$-algebra homomorphism $\Pi_C: A \rightarrow C$ intertwining the comultiplications. One then also has a notion of associated quantum homogeneous space:
\[
B  = \mcO(\K\backslash \U) = \{f\in A \mid (\Pi_C\otimes \id)\Delta_A(f) = f\otimes 1_C\}.
\]
As expected, this $*$-algebra $B$ has a coaction by $A$, i.e.\ $\U$ acts by right translation on $\K\backslash \U$, via
\begin{equation}\label{EqCoactForm}
\alpha_B = \Delta_{\mid B}: B \rightarrow B\odot A,
\end{equation}
where $\odot$ denotes the algebraic tensor product (over $\C$). From classical intuition, one might suspect that any unital \emph{right coideal $*$-subalgebra} $B \subseteq A$, i.e.\ any unital $*$-subalgebra to which $\Delta_A$ restricts to a map as in  \eqref{EqCoactForm}, arises from a compact quantum subgroup. This turns out to be wrong! One can however still make sense of an associated `stabilizer subgroup' $\K$, but one is  forced to go beyond the concept of (quotient) Hopf $*$-algebras (see Section \ref{SecCoideals}).

The above theory of CQG Hopf $*$-algebras shines through the richness of its examples, particularly in the setting of \emph{$q$-deformation of connected, simply connected semisimple Lie groups $U$ and their Lie algebras} $\mfu$ \cite{Dri87,Jim85,LS91}. This richness persists when considering coideal theory and deformations of particular subgroups $K\subseteq U$, namely those $(K,U)$ forming a \emph{symmetric pair}, where $K$ arises as the fixed points of a Lie group involution of $U$. Initial investigations focused on the rank $1$ case $U= SU(2)$ with its subgroup $SO(2)_t$, where $SO(2)_t$ is a conjugate, depending on a real parameter $t$, of the ordinary inclusion $SO(2) \subseteq SU(2)$ \cite{Koo90,NoMi90,Koo93}. Especially the approach through \emph{twisted primitive elements} \cite{Koo93} clarified the tight connection to coideal theory, and led progressively to the development of quantizations of general symmetric pairs, see e.g.\ \cite{NS95,Nou96,Let99,Let02,DCM20} (and see also \cite{Kol14} for the further generalisation into the Kac-Moody setting). 

The above theory of quantization of compact Lie groups can in fact be approached in two ways: either by considering the deformation $\mcO_q(U)$ of the regular function algebra of $U$, which will lead to a CQG Hopf $*$-algebra, or by considering the dual \emph{quantized enveloping algebra} $U_q(\mfu)$, which will only give a Hopf $*$-algebra. The same duality persists in the setting of symmetric pairs, where on the function algebra side one finds as a right coideal $*$-subalgebra a quantization $B = \mcO_q(K\backslash U) \subseteq \mcO_q(U)$ (geometrically, such $K\backslash U$ can be seen as \emph{(pointed) symmetric spaces of compact type}), while on the universal enveloping side one finds a quantization $U_q(\mfk) \subseteq U_q(\mfu)$ as a left coideal $*$-subalgebra, with $\mfk$ the Lie algebra of $K$.

In \cite{DeC24}, it was suggested that a particular \emph{Drinfeld double $*$-algebra}, constructed from the building blocks $\mcO_q(K \backslash U)$ and $U_q(\mfk)$, can serve as a quantization $U_q(\mfl)$ of the unique real Lie algebra $\mfl$ such that 
\begin{itemize}
\item the maximal compact Lie subalgebra of $\mfl$ is isomorphic to $\mfk$, and
\item the complexification of $\mfl$ is isomorphic to the complexification of $\mfu$. 
\end{itemize}
In \cite{DCDT24a}, this was studied in detail for $U =  SU(2)$ and $K = SO(2)_t$, in which case $\mfl = \mfsl(2,\R)_t$ (arising naturally as a unitary conjugation of $\mfsl(2,\R)$ within $\mfsl(2,\C)$). In this case the irreducible unitary representations of $SL_q(2,\R)_t$, interpreted as `admissible' $*$-representations of $U_q(\mfsl(2,\R)_t)$, were classified. 

In this paper, we will consider in general how, for a given coideal $*$-subalgebra $B \subseteq A$ of a CQG Hopf $*$-algebra $A$, one can obtain $*$-representations of its Drinfeld double either by induction of $*$-representations of $B$, or by induction of $*$-representations of its orthogonal dual coideal. Applied in the case of quantum $SL(2,\R)_t$, this will give us a notion of \emph{principal series representations}. By determining the decomposition of the regular representation of quantum $SL(2,\R)_t$, we then obtain that the latter decomposes into principal series as well as discrete series, similar to the case of classical $SL(2,\R)$. 

The precise contents of this paper are as follows. In the \emph{first section}, we recall some general theory on CQG Hopf $*$-algebras and their coideal $*$-subalgebras. In the \emph{second section}, we consider how, both from an algebraic and an analytic perspective, representations of coideals can be induced to representations of their Drinfeld double. In the \emph{third section}, we formally introduce the regular representation, as well as a variant thereof. In the next sections, we then consider the particular case of quantum $SL(2,\R)_t$. In the \emph{fourth section}, we determine how irreducible representations of quantum $SL(2,\R)_t$ decompose into irreducibles for the coideals that it is made up of.  In the \emph{fifth section}, we consider how principal series representations of quantum $SL(2,\R)_t$ can be obtained through induction. Finally, in the \emph{sixth section}, we determine the decomposition of the regular representation of quantum $SL(2,\R)_t$ into irreducibles. We end the paper with a short conclusion and outlook. 

\emph{Conventions}

We denote by $\odot$ the algebraic tensor product (over $\C$). We denote by $\Sigma$ the flip map between tensor products.

When $V,W$ are pre-Hilbert spaces, a linear map $T: V\rightarrow W$ is \emph{adjointable} if there exists a (necessarily unique) linear map $T^*: W\rightarrow V$, referred to as the \emph{adjoint}, with 
\[
\langle w,Tv\rangle = \langle T^*w,v\rangle,\qquad \forall v\in V, w \in W.
\]
When $V,W$ are Hilbert spaces, this is equivalent with $T$ being a bounded operator. We then write in this case $\mcB(V,W)$ for the space of all bounded operators, and $\mcB(V) = \mcB(V,V)$.

If $B$ is an algebra, we write its unit (if it exists) as $1_B$. We write ${}_B\Mod$, resp.\ $\Mod_B$ for the category of left, resp.\ right $B$-modules. If $M \in \Mod_B$ and $N\in {}_B\Mod$, we write $M \odot_B N$ for their balanced tensor product.

Similarly, if $C$ is a coalgebra, we denote by $\varepsilon_C$ its counit. We denote by ${}^C\Mod$, resp.\ $\Mod^C$ the category of left, resp.\ right $C$-comodules. If $(V,\delta_V) \in \Mod^C$ and $(W,\delta_W)\in {}^C\Mod$, then $V\square^CW $ is the \emph{cotensor product}
\[
V \square^C W = \{z \in V\odot W \mid (\delta_V\otimes \id)z = (\id\otimes \delta_W)z\}.
\]
We use standard Sweedler notation for a coalgebra $(C,\Delta_C)$ and a right comodule $(V,\delta_V)$ for it:
\[
\Delta_C(c) = c_{(1)}\otimes c_{(2)},\quad (\Delta_C\otimes \id)\Delta_C(c) = c_{(1)}\otimes c_{(2)}\otimes c_{(3)},\quad \delta_V(v) = v_{(0)}\otimes v_{(1)},\quad \textrm{etc.}
\]

To indicate which structures act on a given vector space $V$, we use bullet notation. E.g.\ ${}_{\bullet}V^{\bullet}$ indicates that there is a left module structure on $V$ and a right comodule structure on $V$, while the notation ${}_{\bullet}V\odot {}_{\bullet}W$ indicates that the relevant module structure is via a natural diagonal action.

If $V$ is a Banach space and $S \subseteq V$ a subset, we write $[S]$ for the norm-closure of the linear span of $S$. 

\section{Preliminaries}\label{SecCoideals}

We refer to \cite{NT14} for a general overview of the (algebraic and analytic) aspects of compact quantum groups, borrowing some specific aspects more directly from \cite{DK94}. We also recall the notion of \emph{Doi-Koppinen datum} and \emph{Doi-Koppinen module}, as introduced in \cite{Doi92, Kop95}, and specifically in the setting of CQG Hopf $*$-algebras in \cite{DeC24}. We refer to the latter paper for details on the further facts on unitary Doi-Koppinen data mentioned in this section.

\subsection{Unitary representations and unitary corepresentations}

In the introduction, we already recalled the notion of a $*$-algebra $B$. It is convenient to use here the following convention on nomenclature: whenever $V$ is a (left) $B$-module, we say that $V$ is a \emph{unitary} $B$-module, or a \emph{unitary $B$-representation}, if moreover $V$ is a pre-Hilbert space on which $B$ acts as adjointable operators satisfying 
\[
\langle v,bw\rangle= \langle b^*v,w\rangle,\qquad \forall v,w\in V,b\in B
\]
for the given $*$-structure on $B$. We call it a \emph{Hilbert $B$-representation} if $V$ is a Hilbert space.  Borrowing the terminology of $*$-category from \cite{Chi18}*{Definition 1.3}, we then denote
\[
{}_B\msRep = \textrm{ the $*$-category of unitary }B\textrm{-representations with adjointable intertwiners as morphisms.}
\] 
If $B$ is non-unital, we assume such modules/representations are algebraically non-degenerate, i.e.\ $BV = V$. We similarly write, using the notion of W$^*$-category as in \cite{GLR85}, 
\[
{}_B\Hilb = \textrm{ the W$^*$-category of Hilbert }B\textrm{-representations with bounded intertwiners as morphisms},
\] 
where for $B$ non-unital we now assume analytic non-degeneracy $[BV] = V$.

It is convenient to consider also the dual notions for coalgebras. Namely, a coalgebra $(C,\Delta_C)$ is called a \emph{$\dag$-coalgebra} if it is equipped with an anti-linear involution $\dag: C \rightarrow C$ such that 
\[
\Delta(c^{\dag}) = \Delta^{\opp}(c)^{\dag\otimes \dag},\qquad c\in C. 
\]
A (right) $C$-comodule $(V,\delta_V)$ is \emph{unitary}, or a unitary corepresentation, if $V$ is a pre-Hilbert space and 
\[
\langle v,w_{(0)}\rangle w_{(1)} = \langle v_{(0)},w\rangle v_{(1)}^{\dag},\qquad \forall v,w\in V. 
\]
We denote by $\Rep^C$ the resulting $*$-category. A $C$-comodule $V$ is called \emph{unitarizable} if it admits a pre-Hilbert space structure for which it becomes unitary. 

E.g., if $(A,\Delta_A)$ is a CQG Hopf $*$-algebra with antipode $S_A$, then $A$ becomes a $\dag$-coalgebra through $\Delta_A$ and
\[
a^{\dag} := S_A(a)^*,\qquad a\in A. 
\] 
The associated notion of unitary $A$-corepresentation is the usual one in the finite-dimensional setting, by associating to $(V,\delta_V)$ the unitary matrix
\[
U\in \mcB(V) \otimes A,\qquad U(v\otimes a) = v_{(0)}\otimes v_{(1)}a,\qquad v\in V,a\in A. 
\]

\subsection{Unitary Doi-Koppinen data}\label{SecUnitaryDK}

A \emph{unitary Doi-Koppinen datum} (= DK-datum) \cite{DeC24}*{Definition 2.12} consists of a triple $(A,B,C)$ with $A$ a CQG Hopf $*$-algebra, $B = (B,\alpha_B)$ a (unital) right $A$-comodule $*$-algebra and $C$ a left $A$-module $\dag$-coalgebra. The compatibility with the $*$-operation and $\dag$-operation are given by 
\[
\alpha_B(b^*) = \alpha_B(b)^{*\otimes *},\qquad (ac)^{\dag} = a^{\dag} c^{\dag},\qquad a\in A,b\in B,c\in C. 
\]
We say that $(A,B,C)$ is of \emph{coideal type} if we are moreover given an equivariant embedding of unital $*$-algebras $B \subseteq A$ and a quotient map of $A$-module $\dag$-coalgebras $\Pi_C: A\twoheadrightarrow C$ such that 
\[
\Ker(\Pi_C) = AB_+,\qquad \textrm{where }B_+ := \Ker(\varepsilon_A)\cap B.
\]
In fact, from any unital right coideal $*$-algebra $B \subseteq A$ one constructs the DK-datum of coideal type $(A,B,C)$ where $C = A/AB_+$ endowed with the quotient $\dag$-coalgebra structure. One can then show that
\begin{equation}\label{EqBAsStab}
B = \{a\in A \mid \Pi_C(a_{(1)})\otimes a_{(2)} = 1_C \otimes a\},
\end{equation}
with $1_C = \Pi_C(1_A)$ the \emph{distinguished grouplike element} of $C$.

\begin{Exa}\label{ExaTrivial}
There are always the following trivial examples of unitary DK-data of coideal type: 
\begin{enumerate}
\item $(A,B,C)$ with $B = A$ (via the identity map) and $C = \C$ (as quotient of $A$ via $\varepsilon_A$). 
\item $(A,B,C)$ with $B = \C1_A$ and $C = A$ (as quotient of $A$ via the identity map). 
\end{enumerate}
\end{Exa}
\begin{Exa}\label{ExaCocomm}
Let $\Gamma$ be a (discrete) group, and $A = \C[\Gamma]$ its group algebra viewed as a CQG Hopf $*$-algebra via 
\[
g^* = g^{-1}, \quad \Delta(g) = g\otimes g,\qquad g\in \Gamma \subseteq \C[\Gamma].
\]
Then if $\Lambda \subseteq \Gamma$ is a subgroup, we obtain $B = \C[\Lambda] \subseteq \C[\Gamma]$ as a right coideal $*$-subalgebra (and in fact Hopf $*$-subalgebra), with associated quotient $C = \C[\Gamma/\Lambda]$ with the natural left $\C[\Gamma]$-module structure and coproduct 
\[
\Delta(g\Lambda) = g\Lambda\otimes g\Lambda,\qquad g\Lambda \in \Gamma/\Lambda \subseteq \C[\Gamma/\Lambda]. 
\]
Any unitary DK-datum of coideal type with $A$ cocommutative can be shown to be of this form.  Note that $C$ is itself a Hopf $*$-algebra (as quotient Hopf $*$-algebra of $A$) if and only if $\Lambda$ is a \emph{normal} subgroup of $\Gamma$. 
\end{Exa}
\begin{Exa}\label{ExaCompactClas}
Let $U$ be a compact Hausdorff group, and let $A = \mcO(G)$ be its associated CQG Hopf $*$-algebra (as in the introduction). If $K \subseteq U$ is a compact subgroup, then $C = \mcO(K)$ arises as a quotient left $A$-module $\dag$-coalgebra (and in fact quotient Hopf $*$-algebra of $A$) via the restriction map 
\[
\pi_C: \mcO(U) \rightarrow \mcO(K),\quad f \mapsto f_{\mid K}. 
\]
One can then identify $B$ (as defined through \eqref{EqBAsStab}) with 
\[
B = \mcO(K\backslash U)  = \{f\in \mcO(U) \mid \forall k\in K,u\in U: f(ku) = f(u)\},
\]
and $(A,B,C)$ becomes a unitary DK-datum of coideal type. In fact, any unitary DK-datum of coideal type with $A$ commutative can be shown to be of this form.  Note that $B$ is itself a Hopf $*$-algebra (as Hopf $*$-subalgebra of $A$) if and only if $K$ is a \emph{normal} subgroup of $U$. 
\end{Exa}
More elaborate examples will be considered in Section \ref{SecYD} and Section \ref{SubSecQSL2R}. 

If $(A,B,C)$ is a unitary DK-datum of coideal type, then $A$ is faithfully flat as a left and right $B$-module. Moreover, $C$ is cosemisimple, and any of its right comodules is unitarizable, uniquely so (up to a scalar) if the comodule is simple. In the following, we always assume fixed a maximal family $\{\msG_{\beta}\}_{\beta\in \Irr(C)}$ of non-equivalent simple = irreducible unitary right $C$-corepresentations, which are then necessarily finite-dimensional. We choose as one of these models the distinguished \emph{trivial} unitary $C$-corepresentation $\msG_{\epsilon} = \C 1_C$, where $1_C$ has norm $1$. See \cite{DeC24}*{Section 2} for more information on these facts.

For $\beta\in \Irr(C)$ and $\xi,\eta\in \msG_{\beta}$, we write the associated matrix coefficient as 
\[
U_{\beta}(\xi,\eta) = (\xi^*\otimes \id)\delta_{\beta}(\eta)\in C, 
\]
and we use the similar notation for general unitary comodules. We write $C_{\beta}$ for the $\dag$-subcoalgebra
\[
C_{\beta} = \mathrm{ span } \{U_{\beta}(\xi,\eta)\mid \xi,\eta\in \msG_{\beta}\}\subseteq B,
\]
and given any unitary $C$-corepresentation $(V,\delta_V)$ we refer to 
\[
V_{\beta} = \{v\in V \mid \delta_V(v) \in V \odot C_{\beta}\}
\] 
as the $\beta$-spectral subspace. We use the similar notation for unitary $A$-comodules, where now $\{\Hsp_{\alpha}\}_{\alpha \in \Irr(A)}$ is a maximal family of non-equivalent irreducible right unitary $A$-corepresentations. 

If $(A,B,C)$ is a unitary DK-datum, we use the notation ${}_B\msMod^C$ for its \emph{Doi-Koppinen modules} = DK-modules, i.e.\ left $B$-modules $V$ endowed with a right $C$-comodule structure $\delta_V$ such that
\[
\delta_V(bv) = \alpha_B(b)\delta_V(v),\qquad \forall b\in B,v\in V. 
\] 
We write ${}_B\msRep^C$ for the associated $*$-category of unitary DK-representations, i.e.\ DK-modules with pre-Hilbert space structure which are at the same time a unitary $B$-representation and unitary $C$-corepresentation.

When we have a unitary DK-datum $(A,B,C)$ and $B_0\subseteq B$ an $A$-comodule $*$-subalgebra and $C\twoheadrightarrow C_0$ an $A$-module quotient $\dag$-coalgebra, then also $(A,B_0,C),(A,B,C_0)$ and $(A,B_0,C_0)$ are unitary Doi-Koppinen data, and we have (joint) restriction and corestriction functors
\[
{}_B\msMod^C \rightarrow {}_{B_0}\msMod^C,\qquad {}_B\msMod^C\rightarrow {}_B \msMod^{C_0},\qquad {}_B\msMod^C \rightarrow {}_{B_0}\msMod^{C_0},
\]
with similar restriction functors for the associated unitary representation categories.

\begin{Exa}
Let us resume Example \ref{ExaCocomm}, and consider a pair of discrete groups $\Lambda \subseteq \Gamma$ with 
\[
A= \C[\Gamma],\quad  B = \C[\Lambda],\quad C = \C[\Gamma/\Lambda].
\]
Then the category ${}_B\msMod^C$ of DK-modules becomes the category of $\Gamma/\Lambda$-graded vector spaces $V = \oplus_{x\in \Gamma/\Lambda}V_x$, endowed with a $\Lambda$-action such that 
\[
g\cdot V_x \subseteq V_{gx},\qquad g\in \Lambda,x\in \Gamma/\Lambda.  
\]
\end{Exa}
\begin{Exa}
Let us resume Example \ref{ExaCompactClas}, and  consider a pair of compact Hausdorff groups $K \subseteq U$ with 
\[
A=  \mcO(U),\quad  B = \mcO(K\backslash U),\quad C = \mcO(K).
\]
Let $g \in U$, and consider $X_g := K \backslash KgK \subseteq K \backslash U$. Then $V = \mcO(X_g)$ may be seen as a quotient $*$-algebra of $\mcO(K \backslash G)$ via restriction, and $\mcO(X_g)$ becomes a DK-module via 
\[
f \cdot h = f_{\mid \mcO(X_g)}h,\qquad (\id\otimes \ev_{k})\delta_V(f) = f(-k),\qquad k \in K,f\in \mcO(K\backslash U),h \in \mcO(X_g).  
\]
\end{Exa}

\subsection{Orthogonal duals of coideal $*$-subalgebras}\label{SecInfDrinfDo}

Let $(A,B,C)$ be a unitary DK-datum of coideal type. We recall the \emph{full and restricted} dual $*$-algebra of $C$ (= \emph{orthogonal} dual $*$-algebra to $B$). If $(A,B,C)$ is a unitary DK-datum, the \emph{full dual $*$-algebra} is given by
\[
\msI = \Lin_{\C}(C,\C),\qquad \omega\chi = (\omega\otimes \chi)\circ \Delta_C,\qquad \omega^*(c) = \overline{\omega(c^{\dag})},\qquad c\in C,\omega,\chi\in \msI. 
\]
We can identify $\msI$ with a direct product of matrix $*$-algebras:
\[
\msI  =  \prod_{\beta\in \Irr(C)} \mcB(\msG_{\beta}),\quad \omega = (\omega_{\beta})_{\beta},\qquad \omega(U_{\beta}(\xi,\eta)) = \langle \xi,\omega_{\beta}\eta\rangle,\qquad \xi,\eta\in \msG_{\beta}.  
\]
We will also occasionally write the pairing between $C$ and $\msI$ as 
\[
\tau(c,\omega) = \omega(c),\qquad \omega \in \msI,c\in C. 
\]
We denote $p_C \in \msI$ the projection onto $\msG_{\epsilon}$, so $(p_C)_{\beta} = \delta_{\beta,\epsilon}$. We write $\Phi_C$ when viewing $p_C$ as a functional:
\[
\tau(p_C,c) = \Phi_C(c),\qquad c\in C. 
\]
We use the similar notation $p_A$ and $\Phi_A$ for $(A,\Delta_A)$ itself. We further write, by abuse of notation,
\[
\Phi_C(a) = \Phi_C(\pi_C(a)),\qquad a \in A. 
\]

Inside $\msI$, one has the $*$-subalgebra $\mcI$, called \emph{restricted dual}:
\begin{equation}\label{EqDecomDirSumMatr}
\mcI = \oplus_{\beta\in \Irr(C)} \mcB(\msG_{\beta}) = \{\Phi_C(a-)\mid a\in A\},
\end{equation}
where the last equality is \cite{DCDT24b}*{Lemma 2.5}. This algebra is non-unital if $C$ is infinite-dimensional.  

We use the similar notation $\msU = \Lin_{\C}(A,\C)$ for the full linear dual of $A$, and $\mcU \subseteq \msU$ for the restricted dual. Note that we can naturally identify 
\[
\msI \subseteq \msU,\qquad \tau(x,a) = \tau(x,\pi_C(a)),\qquad x\in \msI,a\in A.
\]
However, it is important to realize that in general $\mcI \nsubseteq \mcU$. 

From the identification \eqref{EqDecomDirSumMatr}, one concludes that there is an isomorphism of categories
\begin{equation}\label{EqIsoCatsComodMod}
\Mod^C \rightarrow {}_{\mcI}\Mod,\quad xv := \tau(v_{(1)},x)v_{(0)},\qquad x\in \mcI,v\in V. 
\end{equation}
This similarly leads to an isomorphism of $*$-categories 
\begin{equation}\label{EqIsoCatsCoRepRep}
\Rep^C \rightarrow {}_{\mcI}\Rep.
\end{equation}

Finally, we recall the \emph{extended Heisenberg double} $\mcH(A,\msU)$ of $A$ and $\msU$, which is the unital $*$-algebra generated by $A$ and $\msU$ with commutation relations
\begin{equation}\label{EqCommHeis}
\omega a  = a_{(1)} \omega(a_{(2)}-),\qquad a\in A,\omega\in \msU. 
\end{equation}
As a vector space, it is isomorphic to $A \odot \msU$ through the multiplication map. Inside $\mcH(A,\msU)$, we have
\[
\mcH(A,\msI) = \textrm{lin.\ span}\{a\omega \mid a\in A,\omega \in \msI\},\quad \mcH(A,\mcI) = \textrm{lin.\ span}\{a\omega \mid a\in A,\omega \in \mcI\},
\]
\[
\mcD(B,\msI) = \textrm{lin.\ span}\{b\omega \mid b\in B,\omega \in \msI\},\quad \mcD(B,\mcI) = \textrm{lin.\ span}\{b\omega \mid b\in B,\omega \in \mcI\}
\]
as $*$-subalgebras. This distinguishing notation for the latter is explained in \cite{DCDT24a}*{Section 1}, see also below. We refer to $\mcD(B,\mcI)$ as the \emph{Drinfeld double} of $B$ and $\mcI$. 

Through \eqref{EqIsoCatsComodMod} and \eqref{EqIsoCatsCoRepRep}, we then also have isomorphisms 
\begin{equation}\label{EqCorrModCats}
{}_A\!\Mod^C \cong {}_{\mcH(A,\mcI)}\!\Mod,\qquad {}_B\!\Mod^C \cong {}_{\mcD(B,\mcI)}\!\Mod,
\end{equation}
\begin{equation}\label{EqCorrRepCats}
{}_A\!\Rep^C \cong {}_{\mcH(A,\mcI)}\!\Rep,\qquad {}_B\!\Rep^C \cong {}_{\mcD(B,\mcI)}\!\Rep. 
\end{equation}

\begin{Exa}
Let us resume Example \ref{ExaTrivial}: 
\begin{itemize}
\item In case of $(A,B,C) = (A,A,\C)$, one finds $\mcI = \C \varepsilon_A$ and 
\[
\mcD(B,\mcI) \cong A. 
\] 
\item  In case of $(A,B,C) = (A,\C,A)$, one finds $\mcI = \mcU$ and 
\[
\mcD(B,\mcI) \cong \mcU. 
\] 
\end{itemize}
\end{Exa}

\begin{Exa}
Let us resume Example \ref{ExaCocomm}, and consider a pair of discrete groups $\Lambda \subseteq \Gamma$. With $A= \C[\Gamma]$ and $B = \C[\Lambda]$, we can then identify $\msU = \Lin_{\C}(A,\C)$ with the $*$-algebra $\Func(\Gamma,\C)$ of all $\C$-valued functions on $\Gamma$ (with pointwise multiplication and complex conjugation), and similarly we have $\msI = \Func(\Gamma/\Lambda,\C)$ with embedding
\[
\Func(\Gamma/\Lambda,\C) \rightarrow \Func(\Gamma,\C),\quad f \mapsto f\circ \pi_{\Gamma/\Lambda},\quad\textrm{where}\quad \pi_{\Gamma/\Lambda}: \Gamma \rightarrow \Gamma/\Lambda,\quad g \mapsto g\Lambda. 
\]
Then $\mcD(B,\msI)$ is the usual cross product $*$-algebra $\Lambda \ltimes \Func(\Gamma/\Lambda,\C)$, in the sense that 
\[
gf = f(g^{-1}-)g,\qquad g\in \Lambda,f\in \Func(\Gamma/\Lambda,\C).
\]
In this case $\mcI = \Func_{\fin}(\Gamma/\Lambda)$, the $*$-algebra of \emph{finite support functions} on $\Gamma/\Lambda$.
\end{Exa}

Note that the commutation relation \eqref{EqCommHeis} makes sense for any Hopf $*$-algebra $U$ paired with $A$. When $U$ arises as a quantized enveloping algebra, paired with its  associated function algebra $A$ through a Hopf $*$-algebra pairing $\tau: A\times U  \rightarrow \C$, we refer to $\mcH(A,U)$ as the infinitesimal\footnote{We will use the term `infinitesimal' colloquially to refer to anything related to the (quantized) Lie algebra setting, e.g.\ by describing the action of a Lie algebra on the function algebra of its Lie group as being by `infinitesimal translations'.}   Heisenberg double of $A$ and $U$. If then $I \subseteq U$ is a left coideal $*$-subalgebra, we can also make sense of the subalgebras 
\begin{equation}\label{EqMixedHeis}
\mcH(B,I) \subseteq \mcH(A,I) \subseteq \mcH(A,U),
\end{equation}
which we refer to as the \emph{infinitesimal} Heisenberg double of $B$ and $I$, resp.\ of $A$ and $I$. In this case, we can write the commutation relations more symmetrically as 
\begin{equation}\label{EqCommRelMixedHeis}
xf = \tau(f_{(2)},x_{(1)})f_{(1)}x_{(2)},\qquad f\in B,x\in I.
\end{equation}

If moreover $I$ satisfies
\begin{equation}\label{EqInfVanish}
\tau(f_{(1)},x)f_{(2)} = \varepsilon_U(x)f,\qquad f\in B,x\in I,
\end{equation}
it is meaningful to write also $\mcD(B,I) = \mcH(B,I)$, and to call it the \emph{Drinfeld double}. Indeed, in this case we actually have 
\[
\mcD(B,I) \subseteq \mcD(A,U),
\]
where the latter denotes the \emph{Drinfeld double} Hopf $*$-algebra of $A$ and $U$. Moreover, in this case we also have natural unital $*$-homomorphisms
\[
U \rightarrow \msU,\qquad I \rightarrow \msI,\qquad x \mapsto \tau(-,x),
\]
where $\msU$ is the full linear dual of $A$ and $\msI \subseteq \msU$ is the full linear dual of $C = A/AB_+$. In case the pairing $\tau$ is non-degenerate, the embedding $U \rightarrow \msU$ will in particular be faithful, and we then refer to $\msU$ as a completion of $U$ (along $\tau$). Similarly, if $\tau$ is non-degenerate and $B$ is \emph{defined} as the right coideal $*$-algebra satisfying \eqref{EqInfVanish}, we can refer to $\msI$ as a completion of $I$ (as $I$ will be dense in $\msI = \Lin_{\C}(C,\C)$ for the topology of pointwise convergence, see \cite{DCDT24b}*{Remark 2.6}). We then also have an embedding of unital $*$-algebras 
\[
\mcD(B,I) \subseteq \mcD(B,\msI).
\]

\begin{Exa}
Let us resume Example \ref{ExaCompactClas}, with now $U$ a connected compact Lie group with Lie algebra $\mfu$. In this case, we define $U(\mfu)$ as the usual universal enveloping algebra of the complex Lie algebra $\mfg = \mfu \otimes_{\C}\C$, endowed with the unique Hopf $*$-algebra structure such that 
\[
\Delta(X) = X\otimes 1+ 1 \otimes X,\qquad X^* = -X,\qquad X \in \mfu.
\] 
Then we have a natural non-degenerate pairing of Hopf $*$-algebras between $A = \mcO(U)$ and $U(\mfu)$, uniquely determined via
\[
\tau(f,X) = \frac{d}{dt}_{\mid t=0} f(e^{tX}),\qquad X \in \mfu,f\in \mcO(U). 
\]  
If now $K \subseteq U$ is a connected closed subgroup of $U$ with Lie algebra $\mfk \subseteq \mfu$, we have $I = U(\mfk) \subseteq U_q(\mfu)$ as a left coideal $*$-subalgebra (and in fact Hopf $*$-subalgebra). Then $B$, as defined through \eqref{EqInfVanish}, becomes exactly $\mcO(K \backslash U)$. In this case $\mcD(B,I)$ can be identified with the `infinitesimal' cross product $\mcO(K\backslash U) \rtimes \mfk$,  in the sense that 
\[
Xf -fX = D_X(f),\qquad f\in \mcO(K\backslash U),X \in \mfk,
\]
where $D_X(f) \in \mcO(K\backslash U)$ is given by $D_X(f)(Ku) = \frac{d}{dt}_{\mid t=0}f(Kue^{tX})$. 
\end{Exa}

\begin{Rem}
In case of an inclusion of finite groups $\Lambda \subseteq \Gamma$, one can construct both examples of unitary DK-data of coideal type: 
\[
(A,B,C)  = (\C[\Gamma],\C[\Lambda],\C[\Gamma/\Lambda]),\qquad (A',B',C') = (\mcO(\Gamma),\mcO(\Lambda\backslash \Gamma),\mcO(\Lambda)). 
\]
Denoting by $\mcI = \mcO(\Gamma/\Lambda)$ and $\mcI' = \C[\Lambda\backslash \Gamma]$ the respective (full=restricted) dual $*$-algebras to $C$ and $C'$, we then see that we have an isomorphism of $*$-algebras
\[
 \mcD(B,\mcI) = \Lambda \ltimes \mcO(\Gamma/\Lambda) \cong \mcO(\Lambda/\Gamma) \rtimes \Lambda = \mcD(B',\mcI'),
\]
\[
gf \mapsto gS_A(f),\qquad \textrm{where }S_A(f)(\Lambda g) = f(g^{-1}\Lambda),\qquad g\in \Lambda,f\in \mcO(\Gamma/\Lambda). 
\]
\end{Rem}

\subsection{Modular data}\label{SecModData}

Let $A$ be a CQG Hopf $*$-algebra. We recall that the invariant state $\Phi_A: A \rightarrow \C$ is faithful, meaning $\Phi_A(a^*a)= 0$ can only happen for $a=0$.  Through Cauchy-Schwarz, this condition is equivalent with the property that $\Phi_A(ba)=0$ for all $b\in A$ must entail $a=0$.

The couple $(A,\Delta_A)$ admits a (necessarily unique) modular automorphism, meaning an algebra automorphism 
\begin{equation}\label{EqPropModAut}
\sigma_A: A \rightarrow A,\qquad \textrm{satisfying}\quad \Phi_A(ab) = \Phi_A(b\sigma_A(a)),\qquad \forall a,b\in A. 
\end{equation}

The \emph{modular character of $A$} is then the element $\delta_A\in \msU$ defined through 
\[
\delta_A = \varepsilon_A \circ \sigma_A^{-1}. 
\]
It is indeed a character on $A$, and so in particular $\delta_A\circ S_A = \delta_A^{-1}$. With $\delta_A = (\delta_{A,\alpha})_{\alpha}$, we have that each $\delta_{A,\alpha}$ is positive invertible. In particular,  $\delta_A$ admits complex powers $\delta_A^z = (\delta_{A,\alpha}^z)_{\alpha}$ for $z\in \C$. In terms of $\delta_A$, the modular automorphism $\sigma_A$ can be expressed as 
\begin{equation}\label{EqDefModAut}
\sigma_A(a) = \delta_A^{-1/2}(a_{(1)})a_{(2)}\delta_A^{-1/2}(a_{(3)}),\qquad a \in A.
\end{equation}
One has as well that 
\begin{equation}\label{EqPropAntipod}
S^2_A(a) = \delta_A^{-1/2}(a_{(1)})a_{(2)}\delta_A^{1/2}(a_{(3)}),\qquad a \in A. 
\end{equation}

\section{Induction}

We fix a unitary DK-datum $(A,B,C)$ of coideal type. In this section, we show how to induce unitary $B$-representations and unitary $C$-corepresentations into unitary DK-representations. We note that in the setting of locally compact quantum groups, induction of representations was studied (in varying degrees of generality) in e.g.\ \cite{Kus02,Vae05,KKSS12,Ri22}. In the distinct setting of unitary DK-representations, we will be able to follow a quite direct approach, using mostly algebraic techniques.

\subsection{Algebraic induction}

If $V\in {}_B \msMod^A$ and $W \in {}\msMod^C$, one has ${}_{\bullet} V^{\bullet} \odot W^{\bullet}\in {}_B\msMod^C$ via
\[
b\cdot (v\otimes w) = bv\otimes w,\qquad \delta(v\otimes w) = v_{(0)}\otimes w_{(0)}\otimes v_{(1)}w_{(1)}. 
\] 

\begin{Def}
We call ${}_V\!\Ind(W) = {}_{\bullet}V^{\bullet}\odot W^{\bullet}$ the \emph{$V$-induction} of the $C$-comodule $W$. 
\end{Def}
We obtain in this way the induction functor 
\[
{}_V\!\Ind: \Mod^C \rightarrow {}_B\Mod^C.
\]
For example, as we  have $B\in {}_B\msMod^A$ by left $B$-multiplication and the right $A$-comodule structure $\alpha_B$, we can make the $B$-induction 
\[
{}_B\Ind: \Mod^C \rightarrow {}_B\Mod^C.
\]

\begin{Exa}
Let us resume Example \ref{ExaCocomm}, and consider a pair of discrete groups $\Lambda \subseteq \Gamma$ with 
\[
A= \C[\Gamma],\quad  B = \C[\Lambda],\quad C = \C[\Gamma/\Lambda].
\]
For $x\in \Gamma/\Lambda$, let $\C_x$ be the $C$-comodule $\C$ with $\Gamma/\Lambda$-grading concentrated at $x$. Then we can identify 
\[
{}_B\Ind(\C_x) \cong \C[\Lambda]
\]
as a copy of $\C[\Lambda]$ with the usual left $\C[\Lambda]$-action by multiplication, but with a non-trivial $\Gamma/\Lambda$-grading: 
\[
h \in \C[\Lambda]_{hx},\qquad h \in \Lambda.
\]
\end{Exa}

Similarly, if $V \in {}_B\msMod$ and $W\in {}_A\msMod^C$, one has ${}_{\bullet}V \odot {}_{\bullet}W^{\bullet} \in {}_B\msMod^C$ via
\[
b\cdot (v\otimes w) = b_{(0)}v\otimes b_{(1)}w,\qquad \delta(v\otimes w)  = v\otimes w_{(0)}\otimes w_{(1)}. 
\]

\begin{Def}
If $V \in {}_{B}\msMod$, we call $\Ind^W(V) = {}_{\bullet}V\odot {}_{\bullet}W^{\bullet}$ the \emph{$W$-induction} of the $B$-module $V$. 
\end{Def}
Then we get an induction functor
\[
\Ind^W: {}_B\Mod \rightarrow {}_B\Mod^C. 
\]
Unlike in the previous case, there are now at least two natural candidates to consider for $W$. 

\begin{Exa}
We have $C \in {}_A\Mod^C$ through its coproduct and natural $A$-module structure. This leads to the induction functor
\[
\Ind^C: {}_B\Mod \rightarrow {}_B\Mod^C. 
\]
\end{Exa}
\begin{Exa}\label{ExaAModI} 
Let $\mcI$ be the restricted dual of $C$. To make $\mcI$ into an object of ${}_A\Mod^C$, it is by  \eqref{EqCorrModCats} equivalent to make $\mcI$ into an object of ${}_{\Hsp(A,\mcI)}\Mod$. For this, we put 
\[
a\cdot \chi = \chi(S_A^{-1}(a)-),\qquad \omega\cdot \chi = \omega \chi,\qquad a\in A,\omega,\chi\in \mcI. 
\]
It follows from a direct computation that the compatibility \eqref{EqCommHeis} is satisfied.  
\end{Exa}

In general, the objects $C$ and $\mcI$ are not isomorphic in ${}_A\Mod^C$, and not even in ${}_B\Mod^C$.  We will provide an illustrative example in the next section.

\begin{Exa}\label{ExaIndClass}
Let us resume Example \ref{ExaCompactClas}, and  consider a pair of compact Hausdorff groups $K \subseteq U$ with 
\[
A=  \mcO(U),\quad  B = \mcO(K\backslash U),\quad C = \mcO(K).
\]
Let $Kg \in K \backslash U$, and consider $\C_{Kg}$ as the one-dimensional representation of $\mcO(K \backslash U)$ obtained by evaluation at $Kg$. Then we can identify 
\[
\Ind^C(\C_{Kg}) \cong \mcO(K)
\]
as a copy of $\mcO(K)$ as a right $\mcO(K)$-comodule, but with a shifted $\mcO(K \backslash U)$-module structure via 
\[
f\cdot h = f_{Kg}h,\quad f_{Kg}(k) = f(Kgk),\qquad f\in \mcO(K\backslash U),h \in \mcO(K),k\in K. 
\]

Note that in this case, $C$ and $\mcI$ \emph{can} be identified as elements of ${}_A\Mod^C$, via the Fourier transform 
\[
\mcO(K) \rightarrow \mcI,\quad f \mapsto \left(h \mapsto \int_K f(k^{-1})h(k)\rd k\right). 
\]
\end{Exa}

The following shows that whenever induction as above is possible in two ways, one is free to choose which induction one applies. 

\begin{Lem}\label{LemIsoDKmods}
If $V\in {}_B\Mod^A$ and $W \in {}_A\Mod^C$, there is an isomorphism of DK-modules
\[
{}_V\Ind(W) \cong \Ind^W(V).
\]
\end{Lem}
\begin{proof}
The isomorphism is realized through the pair of inverse maps
\begin{multline}\label{EqIsoDiffMod}
\mcX_{V,W}: {}_{\bullet}V^{\bullet} \odot W^{\bullet} \rightarrow {}_{\bullet}V \odot {}_{\bullet}W^{\bullet},\quad v \otimes w \mapsto v_{(0)}\otimes v_{(1)}w,\\ \mcX_{V,W}^{-1}: {}_{\bullet}V \odot {}_{\bullet}W^{\bullet} \rightarrow {}_{\bullet}V^{\bullet} \odot W^{\bullet},\quad v\otimes w \mapsto v_{(0)}\otimes S_A(v_{(1)})w.
\end{multline}
\end{proof}

\begin{Rem}
The map $\mcX_{B,C} \circ \Sigma: C\odot B \rightarrow B\odot C$ is the associated \emph{entwining structure} for $(B^{\opp},C)$ (\cite{BM98}*{Definition 2.1},\cite{Brz99}*{Example 2.2}). 
\end{Rem}

\subsection{Example: Yetter-Drinfeld modules}\label{SecYD}

An illustrative example of a unitary DK-datum is obtained as follows. Fix a CQG Hopf $*$-algebra $(H,\Delta)$. Let $H^{\opp}$ be a copy of  $H$ as a coalgebra, but with the opposite product and the modified $*$-structure 
\[
h^{\star} = S_H^2(h)^*.
\] 
Then $H^{\opp}$ is again a CQG Hopf $*$-algebra. 

Put $A = H \otimes H^{\opp}$ as a tensor product CQG Hopf $*$-algebra. Then $A$ has the quotient module $\dag$-coalgebra $C = (H,\Delta,*)$ through $\Pi_C(h\otimes k^{\opp}) = hk$. The associated right coideal $*$-subalgebra is
\begin{equation}\label{EqIsoCoid}
H \cong B  = (\id\otimes S_H^{-1})\Delta_H^{\opp}(H) \subseteq A
\end{equation}
(see \cite{CMS97} or \cite{DeC24}*{Example 1.4}). Then ${}_B\Mod^C$ is the same as the category of Yetter-Drinfeld $H$-modules ${}_H\mathrm{YD}^H$, i.e.\ left $H$-modules with right $H$-comodule structure such that 
\[
\delta(hv) = h_{(2)}v_{(0)}\otimes h_{(3)}v_{(1)}S_H^{-1}(h_{(1)}). 
\]

Now under the above identifications, we get that $C = H \in {}_H\mathrm{YD}^H$ as follows: 
\[
h\rhd k = h_{(2)}kS_H^{-1}(h_{(1)}),\qquad \delta(k) = \Delta_H(k) = k_{(1)}\otimes k_{(2)},\qquad h,k\in H.
\]

On the other hand, we get that $\mcI$ becomes the restricted dual of $H$, equipped with the YD-structure through left multiplication with $\mcI$ and the left $H$-module structure
\[
h \rhdb \omega = \omega(S_H^{-1}(h_{(2)})- h_{(1)}),\qquad \omega \in \mcI,h\in H. 
\]

\begin{Prop}\label{PropKac}
We have $H$ and $\mcI$ isomorphic in ${}_H\mathrm{YD}^H$ if and only if $(H,\Delta_H)$ is of Kac type, i.e.\ $\Phi_H$ is tracial.
\end{Prop} 
\begin{proof}
The restricted dual $\mcI$ of $H$ is isomorphic to $H$ as a vector space via 
\[
H \rightarrow \mcI,\quad h \mapsto \Phi_H(-S_H(h)). 
\]
Under this isomorphism, the left $\mcI$-module structure on $\mcI$ is simply the one associated to $H$ as a right $H$-comodule under $\Delta_H$. On the other hand, the left 
$H$-module structure gets transported to\footnote{Note that the modular element arises in this setting as a passage from $S_H$ to $S_H^{-1}$, by means of \eqref{EqPropAntipod}.}
\[
h \rhdb k =  \tau(h_{(2)},\delta_H) h_{(3)}kS_H^{-1}(h_{(1)}),\qquad h,k\in H,
\]

If now $(H,\rhd,\Delta_H)$ and $(H,\rhdb,\Delta_H)$ are isomorphic as objects in ${}_H\mathrm{YD}^H$, they are in particular isomorphic as $H$-comodules. But any such isomorphism must be of the form 
\[
\theta_g:H \rightarrow H,\quad h \mapsto \tau(h_{(1)},g)h_{(2)},
\]
for some invertible element $g \in \msI= \Lin_{\C}(H,\C)$. In particular, $\tau(1_H,g)\neq0$, so we assume $\tau(1_H,g)= 1$. 

If now moreover $\theta_g$ intertwines $\rhd$ and $\rhdb$, we have in particular that 
\[
\varepsilon_H(\theta_g(h \rhd 1_H)) = \varepsilon_H(h \rhdb \theta_g(1_H)),\qquad \forall h\in H.
\]
Writing this out, we get
\[
\varepsilon_H(h) = \tau(h,\delta_H), \qquad \forall h\in H,
\]
i.e.\ $\delta_H = \varepsilon_H$. But this is equivalent to the traciality of $\Phi_H$. 
\end{proof}

\subsection{Analytic induction}\label{SecAnInd}

One can easily lift  induction functors to the level of unitary representations. For if $V\in {}_B \msRep^A$ and $W \in {}\msRep^C$, one has 
\[
{}_V\!\Ind(W) = {}_{\bullet}V^{\bullet}\odot W^{\bullet} \in {}_B\Rep^C
\]
as a tensor product pre-Hilbert space, using that $(ac)^{\dag} = a^{\dag}c^{\dag}$ for $a\in A,c\in C$.  

Similarly, if $V \in {}_B\msRep$ and $W\in {}_A\msRep^C$, one has 
\[
\Ind^W(V) = {}_{\bullet}V \odot {}_{\bullet}W^{\bullet} \in {}_B\msRep^C
\]
as a tensor product pre-Hilbert space. 

We then have the following analogue of Lemma \ref{LemIsoDKmods}. 

\begin{Lem}\label{LemIsoDKReps}
If $V\in {}_B\Rep^A$ and $W \in {}_A\Rep^C$, there is a unitary isomorphism of unitary DK-representations
\[
{}_V\Ind(W) \cong \Ind^W(V).
\]
\end{Lem}
\begin{proof}
It is immediately verified that the map $\mcX_{V,W}$ in \eqref{EqIsoDiffMod} is unitary: for $v,v'\in V$ and $w,w'\in W$ we have
\begin{eqnarray*}
\langle v_{(0)}\otimes v_{(1)}w,v_{(0)}'\otimes v_{(1)}'w'\rangle &=& \langle v_{(0)}\otimes w,v_{(0)}'\otimes v_{(1)}^*v_{(1)}'w'\rangle \\
&=& \langle v_{(0)}\otimes w,v_{(0)}'\otimes S_A(v_{(1)}^{\dag})v_{(1)}'w'\rangle \\
&=& \langle v\otimes w,v_{(0)}'\otimes S_A(v_{(1)}')v_{(2)}'w'\rangle\\
&=& \langle v\otimes w,v'\otimes w'\rangle. 
\end{eqnarray*} 
\end{proof} 

We now consider two natural objects to perform induction with.

\begin{Def}
We denote $L^2_0(B)$ for the vector space $B$ endowed with the pre-Hilbert space structure
\[
\langle b,b'\rangle = \Phi_B(b^*b'),\qquad \textrm{where }\Phi_B = (\Phi_A)_{\mid B}. 
\]
\end{Def}

\begin{Prop}
We have $L^2_0(B) \in {}_B\Rep^A$ via left multiplication and $\alpha_B = (\Delta_A)_{\mid B}$.
\end{Prop} 
\begin{proof}
Unitarity of the left $B$-module structure is immediate by definition. The unitarity of the right $A$-comodule structure on $B$ follows from the invariance property of $\Phi_A$.  
\end{proof}

For the induction of $B$-representations, the situation is less straightforward. In general, it is not clear if either $C$ or $\mcI$ can be turned into an object of ${}_A\Rep^C$. There is however a natural object in ${}_A\Rep^C$ to consider: taking any $x = (x_{\beta})_{\beta}\in \msI$ with all $x_{\beta}$ strictly positive, we endow $\mcI$ with the inner product 
\begin{equation}\label{EqInnProdI0}
\langle y,z \rangle_x := \Psi_x(y^*z),\qquad \Psi_x(y) := \sum_{\beta} \Tr_{\beta}(yx),\qquad y,z\in \mcI.
\end{equation}
If we write $L^2_0(\mcI)$ for $\mcI$ with this inner product and its $\mcI$-module structure through left multiplication, the \emph{canonical unitary representation} constructed in \cite{Vae01} will make $L^2_0(\mcI)$ into an object of ${}_A\Rep^C \cong {}_{\mcH(A,\mcI)}\Rep$, which is then independent (up to a canonical unitary intertwiner) of the precise choice of $x$. 

As we do not want to enter the analytic details involved in this construction, we will consider only a special case where this construction simplifies. 

\begin{Def}
We call $\Psi_x: \mcI \rightarrow \C$ \emph{relatively invariant} if there exists $g\in \msU$ such that 
\begin{equation}\label{EqRelInvInt}
\Psi_x(\omega(a-)) = \tau(a,g)\Psi_x(\omega),\qquad \forall \omega \in \mcI,\forall a \in A. 
\end{equation}
\end{Def}
It is clear that $g$ is uniquely determined through $\Psi_x$, and we also refer to \eqref{EqRelInvInt} as \emph{$g$-invariance}. Then $g = (g_{\alpha})_{\alpha}$ has all $g_{\alpha}$ strictly positive, and $a \mapsto\tau(a,g)$ is a character on $A$ \cite{DCDT24b}*{Section 3}. 

\begin{Rem}
We do not know if for any arbitrary unitary DK-datum $(A,B,C)$ of coideal type, there exists a $\Psi_x: \mcI \rightarrow \C$ which is \emph{relatively invariant} with respect to \emph{some} $g\in \msU$. This will however be true in the examples we consider (see \cite{DCDT24b}*{Corollary 3.14} for a sufficient criterion). 
\end{Rem} 

\begin{Prop}\label{PropRepresentUnitI}
Fix $g\in \msU$, and assume $\Psi_x$ is $g$-invariant. Then with respect to the inner product \eqref{EqInnProdI0}, we obtain 
\[
L^2_0(\mcI) \in {}_A\Rep^C \cong {}_{\mcH(A,\mcI)}\Rep
\] 
via left multiplication by $\mcI$ and 
\begin{equation}\label{EqLeftModA}
\pi_g(a)\omega := \tau(a_{(1)},g^{1/2})\omega(S_A^{-1}(a_{(2)})-),\qquad \omega \in \mcI. 
\end{equation}
\end{Prop} 
\begin{proof}
It follows by a direct computation that \eqref{EqLeftModA} defines a left $A$-module structure on $\mcI$, and that together with the left $\mcI$-multiplication it gives an object in ${}_{\mcH(A,\mcI)}\Mod$. To see that then $L^2_0(\mcI)\in {}_{\mcH(A,\mcI)}\Rep$, it is enough to see that \eqref{EqLeftModA} is $*$-preserving. This follows also from a direct computation: 
\begin{eqnarray*}
\langle \pi_g(a^*)\omega,\chi\rangle &=& \overline{\tau(a_{(1)}^*,g^{1/2}})\Psi_x(\omega(S_A^{-1}(a_{(2)}^*)-)^*\chi) \\
&=& \tau(a_{(1)},g^{-1/2})\Psi_x(\omega^*(a_{(2)}-)\chi)\\
&=&  \tau(a_{(1)},g^{-1/2})\Psi_x((\omega^*(\chi(S_A^{-1}(a_{(3)})-)))(a_{(2)}-))\\
&=& \tau(a_{(1)},g^{1/2})\Psi_x(\omega^*(\chi(S_A^{-1}(a_{(2)})-)))\\
&=& \langle \omega, \pi_g(a)\chi\rangle. 
\end{eqnarray*} 
\end{proof}

As an example, let us again consider the situation as in Section \ref{SecYD}. Namely, if $\hat{\psi}$ is the right invariant functional on $\mcI$ viewed as the dual multiplier Hopf $*$-algebra to $(H,\Delta_H)$ \cite{VDae98}, then 
\[
\hat{\psi}(\omega(a-b)) = \delta_H^{-1}(a) \hat{\psi}(\omega)\varepsilon_H(b),\qquad a,b\in H,\omega\in \mcI. 
\]
Hence $\hat{\psi}$ is a $g = \delta_{H}^{-1}\otimes 1$-invariant functional on $\mcI$. Now by (the proof of) \cite{VDae98}*{Proposition 4.9}, the map  
\[
H \rightarrow \mcI,\quad h \mapsto \Phi_H(-h)
\]
is unitary when endowing $\mcI$ with the inner product $\langle \omega,\chi\rangle = \hat{\psi}(\omega^*\chi)$. It follows that we may identify $\mcI \in {}_B\Rep^C$ with $H$ as a right $H$-comodule through $\Delta_H$, and with the left $H$-module strucure 
\[
\pi_g(h)k = \tau(h_{(2)},\delta_H^{1/2})h_{(3)}kS_H^{-1}(h_{(1)}),\qquad h,k\in H. 
\]
As one sees, this DK-module holds the middle ground between the ones on $H$ and $\mcI$ in Section \ref{SecYD}. Moreover, $L^2_0(\mcI)$ will be non-isomorphic to $H$ or $\mcI$ if and only if $(H,\Delta_H)$ is not of Kac type, by a similar proof as for Proposition \ref{PropKac}.

\section{The regular representation}

Throughout this section, we fix again a unitary DK-datum $(A,B,C)$ of coideal type. 

\subsection{The regular representation}

\begin{Def}\label{DefRegularRep}
The regular DK-representation of $(A,B,C)$ is defined as 
\begin{equation}\label{EqPrimRegDK}
{}_{L^2_0(B)}\Ind(L^2_0(\mcI)) \cong \Ind^{L^2_0(\mcI)}(L^2_0(B)).  
\end{equation}
\end{Def}
So, as a pre-Hilbert space this is simply $L_0^2(B) \odot L^2_0(\mcI)$. Note that we only defined $L^2_0(\mcI)$ in detail when $\mcI$ has a relatively invariant weight $\Psi_x$, but the left hand side of \eqref{EqPrimRegDK} can be defined using just the left $\mcI$-module structure on $L^2_0(\mcI)$ with respect to any $\Psi_x$ (here it is more clear directly that any other choice of $x$ will lead to a unitarily equivalent copy as a unitary left $\mcI$-module). 

The motivation to call this the regular representation comes from \cite{DCDT24b}*{Theorem 2.4}: Again assuming that $\Psi_x$ is $g$-invariant, we can endow $\mcD(B,\mcI)$ with the functional 
\begin{equation}\label{EqCanFunctDD}
\varphi_{\mcD}: yb \mapsto \Psi_x(y)\Phi_B(b),\qquad y \in \mcI,b\in B.
\end{equation}
Then $\mcD(B,\mcI)$ becomes a unitary $\mcD(B,\mcI)$-representation under left multiplication and the inner product 
\[
\langle w,z\rangle = \varphi_{\mcD}(w^*z),\qquad w,z\in \mcD(B,\mcI).
\]
This unitary representation is then unitarily equivalent to $L^2_0(B) \odot L^2_0(\mcI)$ as a unitary $\mcD(B,\mcI)$-representation through 
\[
by \mapsto b\otimes y. 
\]

In any case, in general we will then also write the completion of the above regular representation as 
\[
L^2(\mcD(B,\mcI)) = L^2(B) \otimes L^2(\mcI),
\]
and we write the resulting $*$-representation of $\mcD(B,\mcI)$ on it as
\[
\pi_{\reg} : \mcD(B,\mcI) \rightarrow \mcB(L^2(\mcD(B,\mcI))). 
\]

\subsection{A second canonical representation}

Apart from the regular representation defined in Section \ref{DefRegularRep}, there seems to be another `canonical' DK-representation. The proof of the following lemma is immediate:

\begin{Lem}\label{LemSecCan}
The vector space $A$ becomes a unitary DK-representation through the DK-module structure 
\[
b\cdot a = ba\qquad \delta(a) = a_{(1)}\otimes \pi_C(a_{(2)}),\qquad b\in B,a\in A
\]
and the inner product 
\[
\langle a,b\rangle = \Phi_A(a^*b),\qquad a,b\in A.
\]
\end{Lem}
We write $L^2_0(A)$ for $A$ with this inner product and associated unitary DK-representation. We can view it as the natural DK-representation on $A$ for the unitary DK-datum $(A,A,A)$ by restricting to $B$ and co-restricting to $C$. 

The precise relation between $L^2_0(A)$ and the regular representation $L^2_0(B) \odot L^2_0(\mcI)$ seems somewhat mysterious: there is in general no natural (unitary) isomorphism between them, and to the best of the author's knowledge it is not known under which general or natural conditions these unitary DK-representations are isomorphic, even just as modules (see however Proposition \ref{PropRegAndSecDiscr} for a particular case where this can be completely resolved). The only obvious statement that can be made in general, is that any morphism
\[
\theta: C\rightarrow A
\]
of right $C$-comodules leads to a DK-intertwiner
\[
B\odot C \rightarrow A,\quad b\otimes c \mapsto b\theta(c). 
\]
Such morphisms are of course abundant, as $C$ is cosemisimple. But it is not clear in which case they will lead to a (unitary) equivalence. 

We will see however that $L^2_0(A)$ and $L^2_0(B)\odot L^2_0(\mcI)$ arise in a very similar way. To avoid analytic complications, we will consider them simply as DK-modules, and write them as $A$ and $B\odot C$. Indeed, in the latter case we may replace $\mcI$ by $C$ as $\mcI$ and $C$ are isomorphic as $C$-comodules, and then 
\[
B\odot \mcI = \Ind_B(\mcI)\cong \Ind_B(C) = B\odot C.
\] 

To compare $A$ and $B\odot C$, we will consider the natural \emph{tensor product} structure on ${}_B\msMod^C$. For this, we recall by \cite{Tak79}*{Theorem 1} that ${}_B\msMod^C$ is equivalent (as a $\C$-linear category) to ${}_B\msMod^A_B$ via
\[
F: {}_B\msMod^C \rightarrow {}_B\msMod_B^A,\qquad {}_{\bullet}V^{\bullet} \mapsto {}_{\bullet}V\overset{C}{\square} {}_{\bullet}A_{\bullet}^{\bullet},
\]
with adjoint quasi-inverse given by 
\[
G: {}_B\msMod_B^A \rightarrow {}_B\msMod^C,\qquad M \mapsto M/MB_+.
\]
This allows to introduce on ${}_B\msMod^C$ the tensor product defined by 
\[
V\boxtimes W := G({}_{\bullet}F(V)^{\bullet}\otimes_B F(W)_{\bullet}^{\bullet}). 
\]

\begin{Prop}
\begin{enumerate}
\item The regular DK-module $B\odot C$ is isomorphic to $B \boxtimes C$.
\item The  DK-module $A$ is isomorphic to $C\boxtimes B$.
\end{enumerate} 
\end{Prop} 
\begin{proof}
It is enough to show that 
\[
F(A) \cong F(C)\odot_B F(B),\qquad B\odot C \cong G(F(B)\odot_B F(C)).
\]

It is however immediate that $F(C) \cong {}_{\bullet}A_{\bullet}^{\bullet} \in {}_B\msMod_B^A$ with its canonical structure via the mutually inverse maps
\[
A \rightarrow F(C) =  C\overset{C}{\square} A,\quad a \mapsto \pi_C(a_{(1)})\otimes a_{(2)},\qquad C\overset{C}{\square} A \rightarrow A,\quad \sum_i c_i \otimes a_i \mapsto \sum_i \varepsilon_C(c_i)a_i,
\]
while $F(A) \cong {}_{\bullet}A^{\bullet}\odot B_{\bullet}^{\bullet}$ via the mutually inverse maps 
\[
A \odot B \rightarrow A \overset{C}{\square} A,\quad a\otimes b \mapsto a_{(1)}\otimes a_{(2)}b,\qquad 
A \overset{C}{\square} A \rightarrow A \odot B,\quad \sum a_i \otimes a_i' \mapsto \sum_i a_{i(1)} \otimes S_A(a_{i(2)})a_i',
\]
where the inverse has the correct range by equation \eqref{EqBAsStab}. It then also follows immediately that 
\[
{}_{\bullet}B^{\bullet}\odot B_{\bullet}^{\bullet} \cong F(B) = B\overset{C}{\square} A,\qquad a \otimes b \mapsto a_{(1)}\otimes a_{(2)}b. 
\]
Now 
\[
{}_{\bullet}A_{\bullet}^{\bullet} \odot_B ({}_{\bullet}B^{\bullet}\odot B_{\bullet}^{\bullet})  \cong {}_{\bullet}A^{\bullet}\odot B^{\bullet}_{\bullet},\qquad  a\otimes (b\otimes b') \mapsto ab \otimes b'.
\]
The second item of the proposition now follows by the composition of isomorphisms
\[
F(C) \odot_B F(B) \cong {}_{\bullet}A_{\bullet}^{\bullet} \odot_B ({}_{\bullet}B^{\bullet}\odot B_{\bullet}^{\bullet})\cong {}_{\bullet}A^{\bullet}\odot B_{\bullet}^{\bullet} \cong F(A). 
\]
For the first item of the proposition, we note that by the above we have 
\[
F(B) \otimes_B F(C) \cong ({}_{\bullet}B^{\bullet}\odot B_{\bullet}^{\bullet})\odot_B {}_{\bullet}A^{\bullet}_{\bullet}
\] 
and
\[
({}_{\bullet}B^{\bullet}\odot B_{\bullet}^{\bullet})\odot_B {}_{\bullet}A^{\bullet}_{\bullet} \cong {}_{\bullet}B^{\bullet}\odot A_{\bullet}^{\bullet},\qquad (b\otimes b')\otimes a  \mapsto b\otimes b'a. 
\]
The first item of the proposition then follows since 
\[
G( {}_{\bullet}B^{\bullet}\odot A_{\bullet}^{\bullet}) = \frac{{}_{\bullet}B^{\bullet} \odot A^{\bullet}}{{}_{\bullet}B^{\bullet} \odot (AB_+)^{\bullet}} \cong {}_{\bullet} B^{\bullet} \odot (A/AB_+)^{\bullet} = {}_{\bullet} B^{\bullet}\odot C^{\bullet}. 
\]
\end{proof}

\begin{Rem}
In the specific situation of the example treated in Section \ref{SecYD}, we \emph{will} have $A \cong B\odot C$. Indeed, in this case the monoidal category of Doi-Koppinen modules equals the monoidal category of YD-modules by \cite{DeC24}*{Example 1.4}, and the latter is braided. It is not hard to show that the resulting identification in fact produces a unitary intertwiner $L^2_0(A) \cong L^2_0(B)\odot L^2_0(\mcI)$ in this case.
\end{Rem}

Let us also briefly look at the regular representations in the case of Example \ref{ExaCocomm}.

\begin{Exa}
Consider a pair of discrete groups $\Lambda \subseteq \Gamma$ with 
\[
A= \C[\Gamma],\quad  B = \C[\Lambda],\quad C = \C[\Gamma/\Lambda].
\]
Then we can identify $\C[\Gamma/\Lambda]\cong \mcI$ as (right) $C$-comodules through $g\Lambda \mapsto \ev_{g\Lambda}$, and the regular DK-representation becomes 
\[
\C[\Lambda]\odot \C[\Gamma/\Lambda] \subseteq l^2(\Lambda) \otimes l^2(\Gamma/\Lambda),
\]
endowed with the unitary DK-representation structure 
\[
x\cdot (y \otimes z\Lambda) = xy \otimes z\Lambda,\quad y\otimes z\Lambda \in (\C[\Lambda]\odot \C[\Gamma/\Lambda])_{yz\Lambda},\qquad x, y \in \Lambda,z\in \Gamma.
\]
By the map $y \otimes z \Lambda \mapsto y \otimes z\Lambda$, this is unitarily equivalent to the unitary DK-representation on $
\C[\Lambda]\odot \C[\Gamma/\Lambda]$ with
\begin{equation}\label{EqAlternDK}
x\cdot (y \otimes z\Lambda) = xy \otimes xz\Lambda,\quad y\otimes z\Lambda \in (\C[\Lambda]\odot \C[\Gamma/\Lambda])_{z\Lambda},\qquad x, y \in \Lambda,z\in \Gamma.
\end{equation}

On the other hand, the unitary DK-representation as in Lemma \ref{LemSecCan} becomes 
\[
\C[\Gamma]\subseteq l^2(\Gamma) 
\]
endowed with the unitary  DK-representation 
\[
x\cdot z = xz,\quad z \in \C[\Gamma]_{z\Lambda},\qquad x\in \Lambda,z\in \Gamma.
\]
\end{Exa}

In fact, we can say a bit more in the case of the example above.

First, we describe the associated cross product $\Lambda \ltimes \Func_{\fin}(\Gamma/\Lambda)$ a bit more directly. Namely, consider the partition of $\Gamma$ into double cosets,
\[
\Gamma = \sqcup_{i\in \Lambda\backslash \Gamma/\Lambda} \Lambda g_i \Lambda,
\]
for certain chosen representatives $g_i \in \Gamma$. Then $\Func_{\fin}(\Gamma/\Lambda)$ splits also into a direct sum 
\[
\Func_{\fin}(\Gamma/\Lambda) \cong \oplus_{i\in \Lambda\backslash \Gamma/\Lambda} \Func_{\fin}(\Lambda g_i \Lambda/\Lambda),
\]
with each summand stable under the $\Lambda$-action. It follows that we can write also 
\[
\Lambda \ltimes \Func_{\fin}(\Gamma/\Lambda) \cong \oplus_{i\in \Lambda\backslash \Gamma/\Lambda} \Lambda \ltimes \Func_{\fin}(\Lambda g_i \Lambda/\Lambda).
\]
Now the action of $\Lambda$ on $\Lambda g_i \Lambda/\Lambda$ by left translation is transitive, and the $\Lambda$-stabilizer of $g_i \Lambda$ is given by 
\begin{equation}\label{EqStabAct}
K_i := \Lambda \cap g_i \Lambda g_i^{-1} \subseteq \Lambda. 
\end{equation}
By (an easy version of) the Mackey imprimitivity theorem \cite{Mac49}, we then have a Morita equivalence between non-degenerate $\Lambda \ltimes \Func_{\fin}(\Gamma/\Lambda)$-representations and non-degenerate $\oplus_{i\in \Lambda\backslash \Gamma/\Lambda}\C[K_i]$-representations, obtained by 
\begin{equation}\label{EqDecompDoubleCoset}
\Hsp \mapsto \oplus_{i\in  \Lambda\backslash \Gamma/\Lambda} \delta_{g_i \Lambda}\Hsp, 
\end{equation}
with $\delta_{g_i\Lambda}\in \Func_{\fin}(\Gamma/\Lambda)$ the indicator function at the point $g_i\Lambda$.

We can now prove:
\begin{Prop}\label{PropRegAndSecDiscr}
Let $\Lambda \subseteq \Gamma$ be an inclusion of discrete groups, and consider the associated unitary Doi-Koppinen datum 
as in Example \ref{ExaCocomm}:
\[
A= \C[\Gamma],\quad  B = \C[\Lambda],\quad C = \C[\Gamma/\Lambda].
\]
Then the regular DK-representation and the DK-representation in Lemma \ref{LemSecCan} are unitarily equivalent if and only if 
\begin{equation}\label{EqHeckePairMod}
[\Lambda : g\Lambda g^{-1} \cap \Lambda] = [\Lambda: g^{-1}\Lambda g \cap \Lambda],\qquad \forall g\in \Gamma. 
\end{equation}
\end{Prop}
In particular, this is always satisfied for $\Gamma$ finite, but can fail for $\Gamma$ infinite: as a concrete example, one can consider the action of $\Z$ on $\Q$ generated by the automorphism $q\mapsto 2q$, with then 
\[
\Gamma = \Z \ltimes \Q = \{(m,q)\mid m\in \Z,q\in Q\},\quad \Lambda = \{(0,n)\mid n \in \Z\} \subseteq \Gamma,\quad g = (1,0). 
\] 
\begin{proof}
Under the Morita equivalence \eqref{EqDecompDoubleCoset}, we have that the regular DK-representation, in the presentation \eqref{EqAlternDK}, gets sent to
\[
l^2(\Lambda) \otimes l^2(\Gamma/\Lambda) \mapsto \oplus_{i\in \Lambda \backslash \Gamma/\Lambda} l^2(\Lambda),
\]
where $\oplus_{i\in \Lambda\backslash \Gamma/\Lambda} \C[K_i]$ acts on the $i$-th component simply by multiplication on the left with $K_i$ as defined by \eqref{EqStabAct}. 

On the other hand, this Morita equivalence sends $l^2(\Gamma)$ to 
\[
l^2(\Gamma) \mapsto \oplus_{i\in \Lambda \backslash \Gamma/\Lambda} l^2(g_i\Lambda),
\]
where $\oplus_{i\in \Lambda\backslash \Gamma/\Lambda} \C[K_i]$ again acts on the $i$-th component by multiplication on the left with $K_i$ (note that indeed $K_ig_i\Lambda \subseteq g_i\Lambda$). 

If we then put
\[
K_i' =  \Lambda \cap g_i^{-1} \Lambda g_i,\qquad \theta_i: K_i \overset{\cong}{\rightarrow} K_i',\quad k \mapsto g_i^{-1}kg_i,
\]
we see that the regular DK-representation and the DK-representation in Lemma \ref{LemSecCan} are unitarily equivalent  if and only if for each $i \in \Lambda \backslash \Gamma/\Lambda$ we have
\[
l^2(\Lambda) \cong {}_{\theta_i}l^2(\Lambda)
\]
as $K_i$-representations by left multiplication (twisted by $\theta_i$ on the right). Clearly this is the case if and only if 
\[
[\Lambda : g_i\Lambda g_i^{-1} \cap \Lambda] = [\Lambda: g_i^{-1}\Lambda g_i \cap \Lambda].
\]
Since one of the $g_i$ can always be chosen to be a pre-determined $g\in \Gamma$, this proves the proposition. 
\end{proof}

\begin{Rem}
We recall that $\Lambda \subseteq \Gamma$ is a \emph{Hecke pair} precisely when
\[
 [\Lambda : g\Lambda g^{-1} \cap \Lambda] < \infty,\qquad \forall g\in \Gamma. 
\]
In this case (with the small extra assumption of reducedness, i.e.\ $\cap_{x\in \Gamma} x\Lambda x^{-1}= \{e\}$), $\Gamma$ can be densily embedded into a totally disconnected locally compact group $G$, its \emph{Schlichting completion}, and \eqref{EqHeckePairMod} will hold if and only if the modular function of $G$ is trivia. See e.g.\ \cite{KLQ08} for more information, and \cite{AV16} (in particular the proof of \cite{AV16}*{Proposition 2.7.}) for more on categories related to those of unitary DK-representations in this setting.      
\end{Rem}

\section{Branching rules for quantum $SL(2,\R)$}

Throughout the remainder of the article, we fix a parameter $0<q<1$, and we use the following notations for $x\in \R$: 
\[
[x] := \frac{q^x - q^{-x}}{q-q^{-1}}, \qquad \qbra{x} := q^x - q^{-x},\qquad \qdif{x} := q^x+ q^{-x}. 
\]
We also fix a number $a\in \R$, and put 
\begin{equation}\label{EqPart}
t = \qbra{a} = q^a-q^{-a}.
\end{equation}

In the following sections, we look at some of the theory of the previous sections in the case where $A = \mcO_q(SU(2))$ \cite{Wor87a} and $B = \mcO_q(S_t^2)$, one of the \emph{Podle\'{s} spheres} \cite{Pod87}. In this case, one interprets the associated Drinfeld double coideal as a ($t$-dependent) quantization of the convolution algebra of $SL(2,\R)$.

In this first section, we recall the associated `infinitesimal' $*$-algebra $U_q(\mfsl(2,\R)_t)$ and its $KAN$-decomposition 
\[
U_q(\mfsl(2,\R)_t) \cong U_q(\mfso(2)_t) \odot \mcO_q(S_t^2).
\]
We then look at how the (admissible) irreducible unitary $U_q(\mfsl(2,\R)_t)$-representations decompose into irreducible unitary $U_q(\mfso(2)_t)$- and $\mcO_q(S_t^2)$-representations.

\subsection{Quantum $SL(2,\R)$}\label{SubSecQSL2R}

We recall the CQG Hopf $*$-algebra $\mcO_q(SU(2))$ \cite{Wor87a,Wor87b}. 

\begin{Def}
The CQG Hopf $*$-algebra $(\mcO_q(SU(2)),\Delta)$ is universally generated by the entries of a $2$-by-$2$ unitary corepresentation matrix $U = \begin{pmatrix} \alpha &\beta\\ \gamma& \delta\end{pmatrix}$ satisfying
\[
U^* = \begin{pmatrix} \delta & -q^{-1}\beta \\ -q\gamma & \alpha\end{pmatrix}.
\]
\end{Def}
The above conditions force on the generators the universal relations 
\begin{subequations}\label{EqDefRelSUq}
\begin{gather}
\alpha\beta = q\beta \alpha,\quad \alpha \gamma = q\gamma \alpha,\quad \beta\gamma = \gamma \beta,\quad \beta\delta =q\delta\beta,\quad \gamma \delta = q\delta\gamma,\\
\alpha\delta - q\beta\gamma = 1,\qquad \delta \alpha -q^{-1}\gamma \beta = 1,\\
\alpha^* = \delta,\quad \beta^* = -q\gamma,\quad \gamma^* = -q^{-1}\beta,\quad \delta^* = \alpha. 
\end{gather}
\end{subequations}

We now introduce a particular right coideal $*$-subalgebra $\mcO_q(S_t^2) \subseteq \mcO_q(SU(2))$. Recall that we fixed also a parameter $t$ in \eqref{EqPart}.

\begin{Def}
We define the \emph{Podle\'{s} sphere} (at parameter $t$) \cite{Pod87} as the unital $*$-subalgebra 
\[
\mcO_q(S_t^2) = \mcO_q(SO(2)_t\backslash SU(2))\subseteq \mcO_q(SU(2))
\] 
generated by the matrix entries of 
\begin{equation}\label{EqDefEt}
E_t := U^*L_t U \in M_2(\mcO_q(SU(2))),\qquad \textrm{where }L_t = \begin{pmatrix} 0 & i \\ -i & -t\end{pmatrix}.
\end{equation}
\end{Def}
We can then write
\[
E_t = \begin{pmatrix} q^{-1}Z_t & Y_t \\ X_t & -t-qZ_t\end{pmatrix}
\]
where $X_t,Y_t,Z_t$ satisfy the universal relations  
\begin{subequations}
\begin{gather}
\label{EqCommRelPod1}X_tZ_t = q^2 Z_tX_t,\qquad Y_tZ_t = q^{-2}Z_tY_t,\\
\label{EqCommRelPod2} X_tY_t = 1-qtZ_t -q^{2}Z_t^2,\qquad Y_tX_t = 1- q^{-1}tZ_t -q^{-2}Z_t^2,\\
\label{DefStarQHS}X_t^* = Y_t,\quad Z_t^*= Z_t.
\end{gather}
\end{subequations}

By construction, one sees that $\mcO_q(S_t^2) \subseteq \mcO_q(SU(2))$ is a right coideal $*$-subalgebra, so one can consider the associated unitary DK-datum of coideal type as described in Section \ref{SecUnitaryDK}:
\[
(A,B,C) = (\mcO_q(SU(2)),\mcO_q(S_t^2),\mcO_q(SO(2)_t)).
\] 
We further also write the full and restricted (orthogonal) duals of $\mcO_q(SU(2))$ and $\mcO_q(S_t^2)$ as
\[
\msU = \msU_q(SU(2)),\qquad \mcU = \mcU_q(SU(2)),\qquad \msI = \msU_q(SO(2)_t),\qquad \mcI = \mcU_q(SO(2)_t). 
\]
This notation is modelled after Example \ref{ExaCompactClas}, as indeed the above quantizes the example given there for the inclusion of (a conjugated copy of) $SO(2)$ inside $SU(2)$, see \cite{DCDT24a}*{Section 2.2} for more information.

\begin{Def}
We define 
\[
\mcU_q(SL(2,\R)_t) = \mcD(\mcO_q(S_t^2),\mcU_q(SO(2)_t))
\] 
to be the \emph{convolution $*$-algebra of quantum $SL(2,\R)$} (at parameter $t$). 
\end{Def}
As explained in \cite{DCDT24a}, this construction is motivated through \emph{Drinfeld duality}, as one can interpret 
\begin{equation}\label{EqDrinfDuali}
\mcO_q(S_t^2) = U_q((\mfa \oplus \mfn)_t) = \mcU_q((AN)_t),
\end{equation}
where $(\mfa \oplus \mfn)_t$ is a conjugated copy of the real Lie algebra $\mfa \oplus \mfn = \left\{\begin{pmatrix} h & n \\ 0 & -h\end{pmatrix}\mid  h,n\in \R\right\}$, and where the quantization does not really see the difference between the infinitesimal version $U_q((\mfa\oplus \mfn)_t)$ and the analytic version $\mcU_q((AN)_t)$, as through quantization $AN$ has been `discretized'. In any case, one can then interpret bijectivity of the multiplication map 
\[
\mcU_q(SO(2)_t) \odot \mcU_q((AN)_t) \rightarrow \mcU_q(SL(2,\R)_t)
\]
as an analogue of the $KAN$-decomposition of $SL(2,\R)$. 

It is convenient to have an algebraic model $U_q(\mfsl(2,\R)_t)$ for $\mcU_q(SL(2,\R)_t)$, where $U_q(\mfsl(2,\R)_t)$ is then to be seen as the quantized enveloping $*$-algebra associated to quantum $SL(2,\R)$. 

Recall that $U = U_q(\mfsu(2))$ is the universal Hopf $*$-algebra generated by elements $E,F,K^{\pm1}$ with
\[
KE = q^2EK,\quad KF = q^{-2}FK,\quad EF-FE = \frac{K-K^{-1}}{q-q^{-1}},\qquad K^* = K,\quad E^* = FK,\quad F^* = K^{-1}E, 
\]
and coproduct
\[
\Delta(K) = K \otimes K,\qquad \Delta(E) = E\otimes 1 + K\otimes E,\qquad \Delta(F) = F\otimes K^{-1}+ 1\otimes F. 
\]
Inside, we have the \emph{skew primitive element} $B_t$ \cite{Koo93} and its associated left coideal $*$-subalgebra:
\[
B_ t= q^{-1/2}(E-FK)-i(q-q^{-1})^{-1}tK \in U_q(\mfsu(2)), \qquad I = U_q(\mfso(2)_t) = \C[B_t]\subseteq U_q(\mfsu(2)), 
\]
where $\C[B_t]$ denotes the polynomial algebra generated by $B_t$. Note that $iB_t$ is self-adjoint. We will write the skew-primitivity of $B_t$ in the following specific form: 
\begin{equation}\label{EqSkewPrim}
\Delta(iB_t) = iB_t  \otimes 1 + K \otimes (iB_t-[a]),
\end{equation}

We have a unique non-degenerate pairing\footnote{In the setting of Hopf $*$-algebras, a pairing satisfies the usual duality between products and coproducts, as well as the compatibility $\tau(f^*,x) = \overline{\tau(f,S(x)^*)}$ and similarly $\tau(f,x^*) = \overline{\tau(S(f)^*,x)}$.}  of Hopf $*$-algebras between $U_q(\mfsu(2))$ and $\mcO_q(SU(2))$ via
\begin{equation}\label{EqSpinOneHalf}
\tau(U,K) = \begin{pmatrix} q & 0 \\ 0 & q^{-1}\end{pmatrix},\quad \tau(U,E) = \begin{pmatrix} 0 & q^{1/2} \\ 0 & 0 \end{pmatrix},\quad \tau(U,F) = \begin{pmatrix} 0 & 0 \\ q^{-1/2} & 0 \end{pmatrix}.  
\end{equation}
Then recalling that $t= \qbra{a} = q^a-q^{-a}$, we can write 
\begin{eqnarray*}
\mcO_q(S_t^2) &=& \{f \in \mcO_q(SU(2)) \mid \forall x\in U_q(\mfso(2)_t): \tau(x,f_{(1)})f_{(2)} = \varepsilon(x)f\}\\
&=&  \{f \in \mcO_q(SU(2)) \mid \tau(iB_t,f_{(1)})f_{(2)} = [a]f\}.
\end{eqnarray*}
Through the pairing above, we have dense embeddings of $*$-algebras
\[
U_q(\mfsu(2)) \hookrightarrow \msU_q(SU(2)),\qquad U_q(\mfso(2)_t) \hookrightarrow \msU_q(SO(2)_t),
\]
the density referring to the pointwise topology on $\msU_q(SU(2)) = \Lin_{\C}(\mcO_q(SU(2)),\C)$. 

We can now define $U_q(\mfsl(2,\R)_t)$ to be the infinitesimal Drinfeld double $\mcD(\mcO_q(S_t^2),U_q(\mfso(2,\R)_t)$: recalling the notions introduced at the end of Section \ref{SecInfDrinfDo}, this is the universal unital $*$-algebra generated by $\mcO_q(S_t^2)$ and $U_q(\mfso(2)_t)$ with commutation relations 
\[
xf = \tau(f_{(2)},x_{(1)})f_{(1)}x_{(2)},\qquad f\in \mcO_q(S_t^2),x\in U_q(\mfso(2)_t).  
\]
This leads to an embedding of unital $*$-algebras
\[
U_q(\mfsl(2,\R)_t) \hookrightarrow \msU_q(SL(2,\R)_t),\quad f \mapsto f,\quad x \mapsto \tau(-,x),\qquad f\in \mcO_q(S_t^2),x\in U_q(\mfso(2)_t). 
\]

Then we have an isomorphism of $*$-categories
\[
{}_{\mcU_q(SL(2,\R)_t)}\Rep \cong {}_{U_q(\mfsl(2,\R)_t)}\Rep_{\loc.\comp.},
\]
where the right hand side denotes the $*$-category of \emph{locally complete} unitary $U_q(\mfsl(2,\R)_t)$-representations, that is, unitary representations with $iB_t$ diagonalizable, with norm-complete eigenspaces and with its eigenvalues in $\{[a+n]\mid n \in \Z\}$ (see \cite{DeC24}*{Corollary 4.9}). Inside, we then have the full sub-category 
\begin{equation}\label{EqAdmissib}
{}_{U_q(\mfsl(2,\R)_t)}\Rep_{\adm} \subseteq {}_{U_q(\mfsl(2,\R)_t)}\Rep_{\loc.\comp.}
\end{equation}
of \emph{admissible} unitary $U_q(\mfsl(2,\R)_t)$-representations, where moreover the eigenspaces of $iB_t$ are required to be finite-dimensional (\cite{DCDT24a}*{Definition 3.1}).

\subsection{Irreducible unitary representations of $SL(2,\R)_t$}\label{SecIrrRepSL2}

The irreducible admissible unitary representations\footnote{Here irreducibility is understood in the analytic sense, i.e.\ there are no closed invariant subspaces.} of $U_q(\mfsl(2,\R)_t)$ were classified in \cite{DCDT24a}. This classification is achieved by considering
\begin{itemize}
\item the (real) value $\lambda$ taken in an irreducible by the \emph{Casimir element}
\begin{equation}\label{EqCasimir}
\Omega_t := iq^{-1}X_t + (q-q^{-1})iZ_tB_t -iqY_t,
\end{equation}
which is a central, self-adjoint element in $U_q(\mfsl(2,\R)_t)$, and
\item the spectrum $S$ of $iB_t$ in the irreducible representation, which will be multiplicity-free.
\end{itemize}
Together, these completely determine (and distinguish) the irreducible admissible unitary representations up to unitary equivalence. The following cases can arise, where again we write $t= \qbra{a} = q^a-q^{-a}$, where $S$ refers to the eigenvalues of $iB_t$, and where $\lambda$ always denotes the value that the Casimir takes:
\begin{enumerate}
\item $L_{\lambda}^+$: with $s = \lceil a\rceil-a$, these are the cases where $\lambda$ can be any real number in the interval
\[
-\qdif{1-2s}<\lambda <\qdif{1},\qquad \textrm{with  } S_{\lambda}^+ = \{[a+2k] \mid k \in \Z\}.
\]
We refer to these as the \emph{even principal series} if $-2<\lambda<2$, and as the \emph{complementary series} outside this interval. The reason for this terminology will be motivated by Theorem \ref{TheoPrincip}.
\item $L_{\lambda}^-$: with $s = \lceil a-1/2\rceil-(a-1/2)$, these are the cases where $\lambda$ can be any real number in the interval
\[
 -\qdif{1-2s}< \lambda < \qdif{0},\qquad \textrm{with  }
S_{\lambda}^- = \{[a+2k+1] \mid k \in \Z\}.
\]
We refer to these as the \emph{odd principal series}.
\item $D_n^+$: these are the cases arising, for any $n \in \Z_{>0}$, with
\[
\lambda =\qdif{n-1},\qquad S_{\lambda}^+ = \{[a+n+2k]\mid k \in \Z_{\geq0}\}.
\]
We refer to $D_2^+,D_3^+, ...$ as the \emph{positive discrete series} and to $D_1^+$ as the mock positive discrete representation. The latter terminology will be motivated through Theorem \ref{TheoDecompPrinc}.
\item  $E_n^+$: these are the cases arising, for any $n \in \Z_{>0}$, with 
\[
\lambda = -\qdif{2a+\lfloor -2a\rfloor +n -1},\qquad 
S_{\lambda}^+ = \{[a+\lfloor -2a\rfloor +n+2k] \mid k \in \Z_{\geq0}\}.
\]
We refer to  $E_2^+, E_3^+,...$ as the \emph{exceptional positive discrete series}\footnote{Note that we have shifted the indexing with respect to the one in \cite{DCDT24a}, to have better compatibility with the indexing for the discrete series.}, and to  $E_1^+$ as the exceptional positive mock discrete representation.
\item $D_n^-$: these are the cases arising, for any $n \in \Z_{>0}$, with
\[
\lambda = \qdif{n-1},\qquad S_{\lambda}^- = \{[a-n-2k] \mid k \in \Z_{\geq0}\}.
\]
We refer to $D_2^-,D_3^-, ...$, as the \emph{negative discrete series} and to $D_1^-$ as the mock negative discrete representation.
\item $E_n^-$: these are the cases arising for any $n \in \Z_{>0}$, with 
\[
\lambda = -\qdif{-2a+\lfloor 2a\rfloor + n-1},\qquad 
S_{\lambda}^- = \{[a-\lfloor 2a\rfloor -n-2k] \mid k \in \Z_{\geq0}\}.
\]
We refer to  $E_2^-, E_3^-, ...$ as the \emph{exceptional negative discrete series}, and to $E_1^-$ as the exceptional negative mock  discrete representation.
\item $T^+$, the \emph{trivial representation}, which has 
\[
\lambda = \qdif{1}, \qquad S_{\lambda}^+ = \{[a]\},
\]
\item $T^-$, a one-dimensional representation which only arises in case $a \in \frac{1}{2}\Z$, and which has 
\[
\lambda=  -\qdif{1},\qquad S_{\lambda}^- = \{[-a]\}.
\]
\end{enumerate}
We write $\Hsp = \Hsp_{\lambda,S}$ or more concretely $\mcL_{\lambda}^{\pm},\mcD_n^{\pm}$ etc.\ for their Hilbert space completions, which carry the corresponding irreducible unitary representations of $\mcU_q(SL(2,\R)_t)$. We write the corresponding algebraic part as $V = V_{\lambda,S}$. So, more directly, one has $V = \mcU_q(SO(2)_t)\Hsp$.

The infinite representations $V_{\lambda,S}$ have the following more concrete description \cite{DCDT24a}*{Proposition 3.20}: they have an orthogonal (but not necessarily orthonormal) basis $\xi_{[c]}$ for $[c]\in S$, with the generators of $U_q(\mfsl(2,\R)_t)$ acting via 
\[
iB_t\xi_{[c]} = [c]\xi_{[c]}
\] 
and with
\begin{eqnarray}
\nonumber iX_t\xi_{[c]}&=& q^{c+1}\frac{\qdif{a\pm (c+1)}\pm\lambda}{\qdif{c}\qdif{c+1}}\xi_{[c+2]}+ \frac{\qbra{c}\qbra{a}+\qdif{2}\lambda}{\qdif{c-1}\qdif{c+1}}\xi_{[c]} - q^{-c+1}\frac{\qdif{a\mp (c-1)}\mp \lambda}{\qdif{c-1}\qdif{c}}\xi_{[c-2]}, \\
\label{EqFormActXYZ} Z_t\xi_{[c]} &=& -\frac{\qdif{a\pm(c+1)}\pm \lambda}{\qdif{c}\qdif{c+1}}\xi_{[c+2]} + \frac{-\qdif{2}\qbra{a} +\qbra{c}\lambda}{\qdif{c-1}\qdif{c+1}}\xi_{[c]} - \frac{\qdif{a\mp(c-1)}\mp \lambda}{\qdif{c-1}\qdif{c}}\xi_{[c-2]}, \\
\nonumber iY_t\xi_{[c]} &=& q^{-c-1} \frac{\qdif{a\pm (c+1)}\pm \lambda}{\qdif{c}\qdif{c+1}}\xi_{[c+2]} - \frac{\qbra{c}\qbra{a}+\qdif{2}\lambda}{\qdif{c-1}\qdif{c+1}}\xi_{[c]} -q^{c-1} \frac{\qdif{a\mp(c-1)}\mp \lambda}{\qdif{c-1}\qdif{c}}\xi_{[c-2]},
\end{eqnarray} 
where we interpret $\xi_{[c]}=0$ if $[c]\notin S$. The ambiguity of sign is determined as follows: 
\begin{itemize}
\item both choices of $\pm$ give equivalent representations in case of $L_{\lambda}^{\mu}$,
\item we choose $+$ in case of $D_n^-$ or $E_n^+$, 
\item we choose $-$ in case of $D_n^+$ or $E_n^-$. 
\end{itemize}
The one-dimensional unitary representations $\pi_{T^{\pm}}$ are given by the formulas 
\[
\pi_{T^{\pm}}(X_t) = \mp i,\quad \pi_{T^{\pm}}(Y_t) = \pm i,\quad \pi_{T^{\pm}}(Z_t) = 0,\quad \pi_{T^{\pm}}(iB_t) = [\pm a].
\]

\begin{Rem}\label{RemType1}
As one sees from the list of irreducible $\mcU_q(SL(2,\R)_t)$-representations, every of such irreducible $*$-representations contains a compact operator in its image (e.g.\ a projection on an eigenspace of $iB_t$). It follows that (the universal enveloping C$^*$-algebra of) $\mcU_q(SL(2,\R)_t)$ will be type $I$, i.e.\ any of its non-degenerate $*$-representations can be written as a direct integral of irreducible ones \cite{Dix77}*{Theorem 8.6.6}.
\end{Rem} 

\subsection{Branching rules}

We have that $U_q(\mfso(2)_t)$ is a (free) unital polynomial $*$-algebra generated in the single selfadjoint generator $iB_t$. Its completion $\msU_q(SO(2)_t)$ (in the sense of Section \ref{SecInfDrinfDo}) can be identified with 
\begin{equation}\label{EqDecompProdSOq2}
\msU_q(SO(2)_t)  \cong \prod_{[a+n], n \in \Z} \C,\qquad 
U_q(\mfso(2)_t) \rightarrow \msU_q(SO(2)_t),\quad iB_t \mapsto ([c])_{[c]}.
\end{equation}
We then write $p_{[c]}\in \msU_q(\mfso(2)_t)$ for the associated minimal projection at position $[c]$.

We refer to 
\[
\theta_{[c]}: U_q(\mfso(2)_t) \rightarrow \C,\qquad \theta_{[c]}(iB_t) = [c],\qquad c\in a+ \Z
\] 
as the \emph{admissible  unitary representations/characters} of $U_q(\mfso(2)_t)$. We then have 
\[
xp_{[c]} = \theta_{[c]}(x)p_{[c]},\qquad x\in U_q(\mfso(2)_t),c\in a+\Z. 
\] 

Now from the description in Section \ref{SecIrrRepSL2}, it is immediately clear how the irreducible admissible unitary $U_q(\mfsl(2,\R)_t)$-representations decompose into irreducible admissible unitary $U_q(\mfso(2)_t)$-representations: they simply become a direct sum (without multiplicity) of the characters $\theta_{[c]}$ with $[c]\in S$.  

On the other hand, to see how their restriction to $\mcO_q(S_t^2)$ behaves, we first recall that the irreducible unitary representations of the latter are classified as follows. 

\begin{Prop}\label{PropRepPod}
Let $t = \qbra{a}$. Every unitary representation of $\mcO_q(S_t^2)$ on a pre-Hilbert space is bounded. Moreover, any irreducible unitary representation of $\mcO_q(S_t^2)$ on a Hilbert space is one of the following: 
\begin{itemize}
\item The representations $\pi_+$ and $\pi_-$ on $\Hsp_{\pm} = l^2(\N)$ such that, denoting by $e_0,e_1,\ldots$ the standard basis,
\begin{equation}\label{EqExprActZ}
\pi_{\pm}(Z_t)e_p = \pm q^{2p\mp a+1} e_p,
\end{equation}
\begin{equation}\label{EqExprActXt}
\pi_{\pm}(X_t)e_p = \sqrt{(1-q^{2p})(q^{\mp 2a+2p}+1)}e_{p-1}\qquad \pi_{\pm}(Y_t)e_p = \sqrt{(1-q^{2p+2})(q^{\mp 2a+2p+2}+1)}e_{p+1}.
\end{equation}
\item The characters $\chi_z$ for $z\in \C$ with $|z|=1$, given by 
\begin{equation}\label{EqCharactersB}
\chi_z(Z_t) = 0,\qquad \chi_z(X_t) = z,\qquad \chi_z(Y_t) = \overline{z}.
\end{equation}
\end{itemize}
In particular, any unitary representation of $\mcO_q(S_t^2)$ on a Hilbert space can be written as a direct sum of copies of $\pi_+$ and $\pi_-$ and a direct integral (with multiplicities) of the characters $\chi_z$. 
\end{Prop}

\begin{Theorem}\label{TheoBranching}
Let $\Hsp = (\Hsp_{\lambda,S},\pi)$ be an irreducible unitary $\mcU_q(SL(2,\R)_t)$-representation. Denote $(\Hsp,\pi')$ for its restriction to a unitary $\mcO_q(S_t^2)$-representation. Then we obtain the following decompositions: 
\begin{enumerate}
\item $\pi' \cong \pi_+\oplus \pi_-$ in the cases $\mcL_{\lambda}^{\mu}$,
\item $\pi' \cong \pi_+$ in the cases $\mcD_n^-$ and $\mcE_n^+$,
\item $\pi' \cong \pi_-$ in the cases $\mcD_n^+$ and $\mcE_n^-$,
\item $\pi' \cong \chi_{\mp i}$ in the cases $\mcT^{\pm}$. 
\end{enumerate}
\end{Theorem}
\begin{proof}
We first show that $\Ker(\pi(Z_t)) = 0$, unless we are in the case $\mcT^{\pm}$. This will already show that in the infinite case, $\pi'$ must be a direct sum of copies of $\pi_+$ and $\pi_-$. 

Consider the $\xi_{[c]}$ as in \eqref{EqFormActXYZ}, and let $e_{[c]}$ be chosen unit normalisations for them. Then from \eqref{EqFormActXYZ}, we see that each of the generators $iX_t,Z_t$ and $iY_t$ acts on $\Hsp$ via a tri-diagonal transformation of the form 
\begin{equation}\label{EqGenForm}
\sum_{[c]\in S} \lambda_{[c]} e_{[c]} \mapsto \sum_{[c]\in S} \lambda_{[c]}(u_{[c]}e_{[c+2]}  + v_{[c]}e_{[c]}+ w_{[c]} e_{[c-2]}).
\end{equation}
From the concrete description in \eqref{EqFormActXYZ}, we see moreover that (say) the coefficient $w_{[c]}$ can only be non-zero for finitely many values $[c]$. This implies that the space $\Hsp^0$ of zero-vectors for $Z_t$ must be finite-dimensional. 

Now by the commutation relations between the $X_t,Y_t$ and $Z_t$, the eigenspace $\Hsp^0$ is preserved by both $X_t$ and $Y_t$, with $X_t$ acting on it as a unitary having $Y_t$ as its inverse. In particular, we can choose a non-zero $\xi \in \Hsp^0$ and a unimodular complex number $\alpha$ such that 
\[
Z_t\xi = 0,\quad X_t\xi  = \alpha \xi,\qquad Y_t\xi = \alpha^{-1}\xi. 
\]

Choose now a projection $p_{[c]}$ onto a $[c]$-eigenspace of $iB_t$ such that $p_{[c]} \xi \neq 0$. Since $p_{[c]}\xi$ lies in the algebraic part $V\subseteq \Hsp$, we can apply $\Omega_t$ to it. If then $\Omega_t = \lambda$ on $V$, we find by definition of $\Omega_t$ that
\[
(q^c - q^{-c})Z_tp_{[c]} \xi = (\lambda - iq^{-1}X_t +iqY_t) p_{[c]}\xi.
\]
If we take the scalar product with $\xi$ and use that $Z_t\xi= 0$, we find 
\[
(\lambda - iq^{-1}\alpha +iq \alpha^{-1})\|p_{[c]}\xi\|^2 =0,
\]
hence 
\[
\lambda = iq^{-1}\alpha -iq \alpha^{-1}.
\]
However, since $\lambda$ is real, this forces $\alpha = \mp i$ and $\lambda = \pm(q^{-1} +q)$. These are precisely the values of $\lambda$ that appear only for $\mcT^{\pm}$.  

It remains to determine the multiplicities of $\pi_{\pm}$ in the infinite case. 

Let us consider first the cases $\mcL^{\pm}_{\lambda}$. Looking more carefully at the precise form of \eqref{EqGenForm} for $X_t$ given by \eqref{EqFormActXYZ}, we see that $u_{[c]},w_{[c]}\neq0$ for all $[c]\in S$. It follows that there are at most two linearly independent zero vectors for $X_t$ in $\Hsp$. On the other hand, from \eqref{EqFormActXYZ} we see that  $\langle \langle \xi_{[c]}, Z_t\xi_{[c]}\rangle$ equals $\qbra{c}\lambda-\qdif{2}\qbra{a}$ up to a strictly positive number, and hence takes on positive as well as negative values as $[c]$ ranges in $S$. It follows that $\Hsp$ contains exactly one copy of $\pi_+$ and one copy of $\pi_-$. 

Assume now that $\Hsp$ is of the form $\mcD_{m}^+$ for $m\in \Z_{>0}$, with boundary value $[a+m]\in S$ (as seen from the list in Section \eqref{SecIrrRepSL2}). Then in particular we must have $\lambda = \qdif{m-1}$. By the same reasoning as before, it now follows that $\Hsp$ can only contain a one-dimensional space of zero vectors for $X_t$ (as \eqref{EqFormActXYZ} shows that $X_t$ can now be modelled by a tri-diagonal operator on $l^2(\N)$). In this case, we see from \eqref{EqFormActXYZ} that, for some $\kappa_n>0$,
\[
\kappa_n \langle \xi_{[a+m+2n]},Z_t\xi_{[a+m+2n]}\rangle =  \qbra{a+m+2n}\qdif{m-1}-\qdif{2}\qbra{a},\qquad n \in \Z_{\geq 0}.
\]
As this latter value tends to $- \infty$ as $n$ goes to infinity, it follows that  $\pi$ will contain exactly one copy of $\pi_-$ and no copy of $\pi_+$.

The argument in the other cases is completely similar. 
\end{proof}

\subsection{Induction and principal series representations}

In this section, we will examine what is the induction, in the sense of Section \eqref{SecAnInd}, of the $*$-characters $\chi_z$ of $\mcO_q(S_t^2)$ given by \eqref{EqCharactersB}. 

We first recall that under the inclusion 
\[
U_q(\mfsu(2)) \subseteq \mcU_q(SU(2)),
\]
we have 
\[
K \mapsto \delta^{1/2}, 
\]
where $\delta$ is the modular character for $\mcO_q(SU(2))$.  Then $\mcU_q(SO(2)_t)$ has (up to a normalisation) a  unique $K^{-1}$-invariant integral $\psi_t$ such that 
\[
\psi_t(p_{[c]}) = \qdif{c},\qquad c \in a+\Z,
\] 
see \cite{DCDT24b}*{Theorem 3.1}. In the following, we write 
\[
g = K^{-1}.
\] 

If we now write 
\[
L^2_{q,0}(SO(2)_t) = L^2_0(\mcU_q(SO(2)_t)),
\]
we can identify $L^2_{q,0}(SO(2)_t)$ with $\C[S_a]$ for $S_a= \{[a+m]\mid m\in \Z\}$, where  the orthonormal unit basis vectors $e_{[c]} \in \C[S_a]$ correspond to normalisations of the projections $p_{[c]}$:  
\[
e_{[c]} = \qdif{c}^{-1/2}p_{[c]}.
\]

To identify the induced characters of $\mcO_q(S_t^2)$, we start with some lemmas. We will write the spin $n/2$-representation of $U_q(\mfsu(2))$ as $V_{n/2}$. We identify $V_{1/2}$ concretely with the representation as given by \eqref{EqSpinOneHalf}, with associated basis $e_+,e_-$, so 
\[
Ke_{\pm} = q^{\pm 1}e_{\pm},\qquad E e_+ = 0,\quad Ee_-= q^{1/2}e_+,\qquad Fe_+ = q^{-1/2}e_-,\quad Fe_-=0. 
\] 
Then $V_{n/2}$ can be identified with the $U_q(\mfsu(2))$-representation spanned by $e_+^{\otimes n}$ in $V_{1/2}^{\otimes n}$. 

\begin{Lem}\label{LemTensProdEigVectiB}
The vectors 
\[
\xi_{1/2,\pm}^{[a]} = \qdif{a}^{-1/2}\begin{pmatrix} q^{\pm a/2} \\\mp iq^{\mp a/2} \end{pmatrix}
\]
are unit eigenvectors of $iB_t$ at eigenvalues $[a\pm 1]$ in $V_{1/2}$. 

More generally, the vectors 
\begin{equation}\label{EqEigenvectiBt}
\xi_{n/2,\pm}^{[a]} = \xi_{1/2,\pm}^{[a\pm(n-1)]} \otimes \xi_{1/2,\pm}^{[a\pm(n-2)]} \otimes \ldots \otimes \xi_{1/2,\pm}^{[a\pm 1]} \otimes \xi_{1/2,\pm}^{[a]}
\end{equation}
are unit eigenvectors of $iB_t$ at eigenvalues $[a\pm n]$ in $V_{n/2}$.
\end{Lem}
\begin{proof}
The proof for $V_{1/2}$ is by a direct computation, using that $\pi_{1/2}(iB_t) = \begin{pmatrix} q[a] & i \\ -i & q^{-1}[a]\end{pmatrix}$. The proof that \eqref{EqEigenvectiBt} is an eigenvector in the general case follows from the skew-primitivity \eqref{EqSkewPrim} of $iB_t$ and induction, for if $\xi$ is an eigenvector for $iB_t$ at eigenvalue $[c]$, we see that 
\[
\Delta(iB_{\qbra{a}})(\xi_{1/2,\pm}^{[c]}\otimes \xi) \underset{\eqref{EqSkewPrim}}{=} (iB_{\qbra{a}} +  ([c]-[a])K)\xi_{1/2,\pm}^{[c]}\otimes \xi = iB_{\qbra{c}}\xi_{1/2,\pm}^{[c]} \otimes \xi = [c\pm 1]\xi_{1/2,\pm}^{[c]} \otimes \xi.
\]
Then $\xi_{n/2,\pm}^{[a]}\in V_{n/2}$ follows since no $V_{k/2}$ with $k <n$ contains an eigenvector for $iB_t$ at eigenvalue $[a\pm n]$ \cite{Koo93}*{Theorem 4.3} (see also 
\cite{DCDT24a}*{Theorem 2.1}). 
\end{proof}

The next lemma describes explicitly the unitary $\mcO_q(SU(2))$-representation on $L^2_{q,0}(SO(2)_t)$ as given by Proposition \ref{PropRepresentUnitI}.

\begin{Lem}\label{LemInductionComp}
Put 
\[
f_{[c]} = \qdif{c}^{-1/2}e_{[c]} =  \qdif{c}^{-1}p_{[c]},\qquad c \in a+\Z.
\]
Then for $n \in \Z_{\geq0}$, we can write 
\[
\qdif{a\pm n} \pi_g(\alpha) f_{[a\pm n]} =  q^{\frac{1}{2}+n\pm a} f_{[a\pm (n+1)]} + q^{\frac{1}{2}-n\mp a}f_{[a\pm (n-1)]},
\]
\[
\qdif{a\pm n} \pi_g(\beta) f_{[a\pm n]} =  \pm iq^{\frac{1}{2}}f_{[a\pm (n+1)]} \mp  iq^{\frac{1}{2}}f_{[a\pm (n-1)]},
\]
\[
\qdif{a\pm n} \pi_g(\gamma) f_{[a\pm n]} =  \mp iq^{-\frac{1}{2}}f_{[a\pm (n+1)]} \pm  iq^{-\frac{1}{2}}f_{[a\pm (n-1)]},
\]
\[
\qdif{a\pm n}\pi_g(\delta) f_{[a\pm n]} =  q^{-\frac{1}{2}-n\mp a}f_{[a\pm (n+1)]} + q^{-\frac{1}{2}+n\pm a}f_{[a\pm (n-1)]},
\]
\end{Lem}
In the course of the proof, the symbols $\epsilon,\nu,\kappa,\ldots$ denote elements of $\{+,-\}$. 
\begin{proof}
Write $U = (U_{\epsilon,\nu})_{\epsilon,\nu}$ for the matrix coefficients in $\mcO_q(SU(2))$ associated to $\pi_{1/2}$ with respect to its basis $e_{\pm}$. Using \eqref{EqSkewPrim}, we obtain from (the infinitesimal analogue of) the commutation relation \eqref{EqCommHeis} that, inside $\mcH(\mcO_q(SU(2)),U_q(\mfso(2)_t))$,
\[
iB_t U_{\epsilon,\nu} = \sum_{\kappa} \langle e_{\kappa},iB_te_{\nu}\rangle U_{\epsilon,\kappa} + \sum_{\kappa} \langle e_{\kappa},Ke_{\nu}\rangle U_{\epsilon,\kappa}(iB_t-[a]).
\]
By a similar computation as in the proof of Lemma \ref{LemTensProdEigVectiB}, we then see that, for $c\in a+\Z$, 
\begin{equation}\label{EqCommRelHeiss}
iB_t \pi_g(U_{\epsilon,\nu}) e_{[c]} =  \sum_{\kappa} \langle e_{\kappa},iB_{\qbra{c}}e_{\nu}\rangle  \pi_g(U_{\epsilon,\kappa})e_{[c]}.
\end{equation}
Since $\pi_{1/2}(iB_{\qbra{c}})$ has eigenvalues $\qbra{c\pm1}$, we see that 
\[
\pi_g(U_{\epsilon,\nu})e_{[c]} = \theta_{\epsilon,\nu}^{[c],+} e_{[c+1]} +  \theta_{\epsilon,\nu}^{[c],-}e_{[c-1]}
\]
for certain scalars $\theta_{\epsilon,\nu}^{[c],\kappa} \in \C$. 

Now the above implies as well that, upon viewing $p_{[c]}$ as functionals on $\mcO_q(SU(2))$, 
\[
p_{[c]}(U_{\epsilon,\nu}-) = \kappa_{\epsilon,\nu}^{[c],+} p_{[c+1]} +  \kappa_{\epsilon,\nu}^{[c],-}p_{[c-1]}
\]
for certain scalars $\kappa_{\epsilon,\nu}^{[c],\kappa} \in \C$. Considering the case $\epsilon = \nu = +$ and $[c] = [a\pm n]$ for $n \geq 0$, evaluating in $\alpha^{n-1}$ and realizing that the $iB_t$-eigenvalue $[a\pm n]$ does not appear in the spin $k/2$-representations of $U_q(\mfsu(2))$ for $k<n$, we see that 
\[
p_{[a\pm n]}(\alpha^n) = \kappa_{+,+}^{[a\pm n],\mp}p_{[a\pm (n-1)]}(\alpha^{n-1}). 
\]
From \eqref{EqEigenvectiBt}, we get 
\[
 \kappa_{+,+}^{[a\pm n],\mp} = |\langle e_+,\xi_{1/2,\pm}^{[a+ n-1]}\rangle |^2 = \frac{q^{n-1\pm a}}{\qdif{a\pm (n-1)}}.
\]
Performing a similar computation with respect to $\delta$, we find 
\[
 \kappa_{-,-}^{[a\pm n],\mp} = |\langle e_-,\xi_{1/2,\pm}^{[a+ n-1]}\rangle |^2 = \frac{q^{-n+1 \mp a}}{\qdif{a\pm(n-1)}}.
\]
Since $S^{-1}(\alpha) = \delta$ and $S^{-1}(\delta) = \alpha$, we easily deduce from this that 
\[
\theta_{\epsilon,\epsilon}^{[a\pm n],\mp} = q^{-\epsilon/2} \frac{\qdif{a\pm (n-1)}^{1/2}}{\qdif{a\pm n}^{1/2}}\kappa_{-\epsilon,-\epsilon}^{[a\pm n],\mp} = \frac{q^{\epsilon/2}q^{\mp \epsilon a}q^{-\epsilon n}}{\qdif{a\pm n}^{1/2}\qdif{a\pm(n-1)}^{1/2}}
\]
Now the fact that $\pi_g$ is unitary gives us $\theta_{\epsilon,\epsilon}^{[a \pm n],\pm} = \theta_{-\epsilon,-\epsilon}^{[a\pm (n+1)],\mp}$. We then straightforwardly deduce the formulas for $\pi_g(\alpha)$ and $\pi_g(\delta)$. 

The formulas for $\pi_g(\beta),\pi_g(\gamma)$ are now immediate from \eqref{EqCommRelHeiss}.
\end{proof}

\begin{Theorem}\label{TheoPrincip}
Assume $z= e^{-i \theta}$. Then there is a unitary equivalence of $\mcU_q(SL(2,\R)_t)$-representations
\begin{equation}\label{EqDirectDecomp}
\Ind^{L^2_{q,0}(SO(2)_t)}(\chi_z)\cong \mcL_{2\sin(\theta)}^+ \oplus \mcL_{2\sin(\theta)}^-.
\end{equation}
\end{Theorem} 
\begin{proof}
From the construction presented in Section \ref{SecAnInd}, we see that $\Ind^{L^2_{q,0}(SO(2)_t)}(\chi_z)$ can be identified, as a $\mcU_q(SL(2,\R)_t)$-module, with $\mcU_q(SO(2)_t) \cong \C[S]$, where $S = \{[a+m]\mid m\in \Z\}$, endowed with the left $\mcU_q(SO(2)_t)$-module structure through multiplication, and with $\mcO_q(S_t^2)$-module structure given through
\[
b\cdot \omega = \chi_z(b_{(1)})\pi_g(b_{(2)})\omega,\qquad 
 b\in \mcO_q(S_t^2),\omega \in \mcU_q(SO(2)_t).
\]
In particular, it is directly clear from this description that $\Ind^{L^2_{q,0}(SO(2)_t)}(\chi_z)$ is admissible (in the sense of \eqref{EqAdmissib}).

Now using $E_t\in M_2(\C)\otimes \mcO_q(S_t^2)$ as in \eqref{EqDefEt}, we have 
\[
(\id\otimes \Delta)E_t = U_{13}^*E_{t,12}U_{13},
\]
so 
\[
(\id\otimes \chi_z\otimes \Delta)E_t = E_{t}(z),
\]
where
\[
E_t(z) = U^*L_t(z)U,\qquad L_t(z) = \begin{pmatrix} 0 & \overline{z} \\ z & -t\end{pmatrix}. 
\]
If we denote $\kappa_z(f) = (\chi_z\otimes \id)\Delta(f)$ for $f\in \mcO_q(S_t^2)$, we see that 
\begin{eqnarray*}
\kappa_z(X_t) &=& -q\overline{z} \gamma^2 + z\alpha^2 -t\alpha\gamma,\\
\kappa_z(Z_t) &=& q\overline{z} \delta \gamma -z\beta \alpha + t \beta\gamma\\
\kappa_z(Y_t) &=& \overline{z}\delta^2 - q^{-1}z \beta^2 + q^{-1}t \beta \delta. 
\end{eqnarray*}
Hence extending this to an inclusion
\[
\kappa_z: U_q(\mfsl(2,\R)_t) \rightarrow \mcH(\mcO_q(SU(2)),U_q(\mfso(2)_t)),\quad f  \mapsto \kappa_z(f),\quad x \mapsto x,\qquad f\in \mcO_q(S_t^2),x\in U_q(\mfso(2)_t),
\]
we find that, with $\Omega_t$ the Casimir element of $U_q(\mfsl(2,\R)_t)$ from \eqref{EqCasimir}, 
\[
\kappa_z(\Omega_t)=   -i\overline{z} \gamma^2 + izq^{-1}\alpha^2 -itq^{-1}\alpha\gamma -i\overline{z}q\delta^2 + iz \beta^2 - it \beta \delta +(q-q^{-1})(q\overline{z} \delta \gamma -z\beta \alpha + t \beta\gamma)iB_t.
\]
Then $\kappa_z(\Omega_t)$ will still commute with $iB_t$, hence $\pi_g(\kappa_z(\Omega_t))e_{[c]}$ will be a scalar multiple of $e_{[c]}$. 

It is then easy (if a bit tedious) to compute that, with the $f_{[c]}$ as in Lemma \ref{LemInductionComp},
\[
\pi_g(\kappa_z(\Omega_t))f_{[a\pm n]} = i(z-\overline{z})f_{[a\pm n]}. 
\]
The theorem now follows from the classification of irreducible admissible unitary $U_q(\mfsl(2,\R)_t)$-representations.
\end{proof}

\begin{Rem}
If we write\footnote{We only put in the extra parameters in the notation of the left hand side for consistency.}
\[
U_q(\mfa_t)  =  \C[X,Y\mid X^* = Y, XY =YX = 1],
\]
then under the interpretation $\mcO_q(S_t^2) = U_q((\mfa \oplus n)_t)$ of \eqref{EqDrinfDuali}, one can view the quotient map 
\[
U_q((\mfa\oplus \mfn)_t) \rightarrow U_q(\mfa_t),\qquad X \mapsto X,\quad Y \mapsto Y,\quad Z \mapsto 0
\]
as the analogue of the quotient map 
\[
U(\mfa\oplus\mfn) \rightarrow U(\mfa). 
\]
This explains why Theorem \ref{TheoPrincip} can really be seen as the correct quantum analogue of principal series induction in the classical setting \cite{Kna01}*{Chapter VII}, except that we induced from the analogue of $AN$ in stead of $MAN$, with $M = \{\pm 1\}$ the centralizer of $\mfa$ in $SU(2)$. However, in our setting the centrality of $M$ allows to view this extra piece of information as simply picking out the apropriate parity for the spectrum of $iB_t$, which corresponds (in complete analogy with the classical setting) to picking out the appropriate component of \eqref{EqDirectDecomp}. It remains to explore this construction more conceptually in the setting of quantizations of higher rank real semisimple Lie groups \cite{DeC24} (but see \cite{VY20} and references therein for the complex setting).
\end{Rem} 

\section{Decomposition of the regular representation for quantum $SL(2,\R)$}

In this section, we compute the direct integral decomposition of the regular representation of quantum $SL(2,\R)$. We start with some preparations on the regular representation of $\mcO_q(S_t^2)$, and present a certain localisation of the algebra $\mcO_q(S_t^2)$ that will be convenient for computations. 

\subsection{Regular representation of $\mcO_q(S_t^2)$}

Recall from Proposition \ref{PropRepPod} the Hilbert representations $\pi_{\pm}$ of $\mcO_q(S_t^2)$ on $\Hsp_{\pm} = l^2(\N)$. Putting $\pi= \pi_{+}\oplus \pi_{-}$, the restriction $\Phi$ of the Haar state of $\mcO_q(SU(2))$ to $\mcO_q(S_t^2)$ is then given by 
\[
\Phi(b) = \frac{q^{-1}-q}{\qdif{a}}\Tr(\pi(b)|Z_t|),\qquad|Z_t| \xi_p^{\pm} = q^{2p\mp a+1}\xi_p^{\pm},
\]
see \cite[Proposition 2.4]{DeC12}. We will use as a shorthand 
\[
T = \frac{q^{-1}-q}{\qdif{a}}|Z_t| \in \mcB(\mcH_+ \oplus \mcH_-),\qquad \textrm{so } \Phi = \Tr(\pi(-)T).
\]

In particular, the GNS-space $L^2_q(S_t^2) := L^2(\mcO_q(S_t^2))$ of $\Phi$ can be identified as a $\mcO_q(S_t^2)$-representation with 
\begin{equation}\label{EqIdentGNS}
L^2_q(S_t^2) \cong (\Hsp_+\otimes \Hsp_+) \oplus (\Hsp_-\otimes \Hsp_-), 
\end{equation}
in a unique way such that
\[
1 \mapsto  \xi_T = \frac{(q^{-1}-q)^{1/2}}{\qdif{a}^{1/2}}\sum_{\pm} \sum_{p\in \Z_{\geq0}} q^{\frac{2p \mp a+ 1}{2}} \xi_p^{\pm}\otimes \xi_p^{\pm}. 
\]
Under this identification,  $\mcO_q(S_t^2)$ acts on the right hand side of \eqref{EqIdentGNS} through $\pi_{\pm}$ on the first components. We can below then also view
\[
L^2_{q,0}(S_{t}^2):= L^2_0(\mcO_q(S_t^2)) \subseteq (\Hsp_+\otimes \Hsp_+)\oplus (\Hsp_- \otimes \Hsp_-). 
\] 

We will further write
\begin{equation}\label{EqPreHilb1}
V_{\pm}  = \textrm{ the linear span of the basis vectors }\xi_{p}^{\pm},\qquad V_{\pm}  \subseteq  \Hsp_{\pm},
\end{equation}
and we write 
\begin{equation}\label{EqPreHilb2}
W= W_+ \oplus W_-,\qquad W_{\pm} = V_{\pm}\otimes V_{\pm}.
\end{equation}
Then $W$ is a pre-Hilbert space, dense in $L^2_q(S_t^2)$, and hence non-degenerately paired with $L^2_{q,0}(S_t^2)\subseteq L^2_q(S_t^2)$ through the scalar product $\langle \cdot,\cdot\rangle$ of $L^2_q(S_t^2)$. We write the natural orthonormal basis of $W_{\pm}$ as 
\begin{equation}\label{EqDoubleInd}
\xi_{k,l}^{\pm} = \xi_k^{\pm}\otimes \xi_l^{\pm},\qquad k,l\in \Z_{\geq 0}. 
\end{equation}

\subsection{A localisation of $\mcO_q(S_t^2)$}

In this section, we construct a particular localisation of $\mcO_q(S_t^2)$, realized as a quotient of a deformed quantized enveloping algebra introduced in the next definition (see also \cite{DeC13}). The consideration of this localisation will be important for the results on the regular representation in subsection \ref{RegularRep}, for the localisation will allow for some particular algebraic properties which would otherwise be absent.

\begin{Def}
For $\mu,\nu \in \{\pm\}$, we define $U_q^{\mu,\nu}(\mfsu(2))$ as the unital $*$-algebra generated by elements $K^{\pm1},E,F$ with 
\begin{equation}\label{EqUnivRelUqVar}
KE = q^2EK,\qquad KF = q^{-2}FK,\qquad EF - FE = \frac{\nu K-\mu K^{-1}}{q-q^{-1}}
\end{equation}
and the $*$-structure 
\[
E^* =FK, \qquad F^* = K^{-1}E,\qquad K^* = K.
\]
\end{Def}

The following four lemmas follow by direct computation.
\begin{Lem}
For any $\mu,\nu,\kappa \in \{\pm\}$, there exists a $*$-homomorphism 
\[
\Delta_{\mu,\nu}^{\kappa}: U_q^{\mu,\nu}(\mfsu(2)) \rightarrow U_q^{\mu,\kappa}(\mfsu(2))\otimes U_q^{\kappa,\nu}(\mfsu(2))
\]
such that 
\[
\Delta_{\mu,\nu}^{\kappa}(E) = E\otimes 1 + K\otimes E,\qquad \Delta_{\mu,\nu}^{\kappa}(F) = F\otimes K^{-1} + 1\otimes F,\qquad \Delta_{\mu,\nu}^{\kappa}(K) = K\otimes K,
\]
satisfying the obvious coassociativity condition. 
\end{Lem}

Note that $U_q^{++}(\mfsu(2)) = U_q(\mfsu(2))$. In the following, we will use the following Sweedler notation: 
\[
\Delta_{\mu,\nu}^{\kappa}(x) = x_{(\mu,\kappa;1)}\otimes x_{(\kappa,\nu;2)},\qquad x\in U_q^{\mu,\nu}(\mfsu(2)).
\]

\begin{Lem}
There exist $*$-homomorphisms $\varepsilon_{\mu}: U_q^{\mu,\mu}(\mfsu(2)) \rightarrow \C$ such that for all $\mu,\nu$ and $x\in U_q^{\mu,\nu}(\mfsu(2))$: 
\[
\varepsilon_{\mu}(x_{(\mu,\mu;1)})x_{(\mu,\nu;2)} = x= x_{(\mu,\nu;1)}\varepsilon_{\nu}(x_{(\nu,\nu;2)}).
\]
Concretely, we have 
\[
\varepsilon_{\mu}(E) = \varepsilon_{\mu}(F) = 0,\qquad \varepsilon_{\mu}(K) = 1.
\]
\end{Lem} 

\begin{Lem}
There exist bijective anti-algebra homomorphisms 
\[
S_{\mu,\nu}: U_q^{\mu,\nu}(\mfsu(2)) \rightarrow U_q^{\nu,\mu}(\mfsu(2))
\]
such that the following is satisfied for all $\mu,\nu$ and $x\in U_q^{\mu,\mu}(\mfsu(2))$: 
\[
S_{\mu,\nu}(x_{(\mu,\nu;1)})x_{(\nu,\mu;2)} = \varepsilon_{\mu}(x)1 = x_{(\mu,\nu;1)}S_{\nu,\mu}(x_{(\nu,\mu;2)}).
\]
Concretely, we have 
\[
S_{\mu,\nu}(E) = -K^{-1}E,\quad S_{\mu,\nu}(F) = -FK,\quad S_{\mu,\nu}(K) = K^{-1}. 
\]
\end{Lem} 

\begin{Lem}
The $*$-algebra $U_q^{+,-}(\mfsu(2))$ becomes a left $U_q(\mfsu(2))$-module $*$-algebra\footnote{This means we have a module algebra for which the $*$-structure satisfies $(y\rhd x)^* = S(y)^*\rhd x^*$.} by the \emph{adjoint action}
\[
y \rhd x =y_{(+,-;1)}x S_{-,+}(y_{(-,+;2)}).
\]
\end{Lem} 

Let us now relate these $*$-algebras to $\mcO_q(S_t^2)$.

\begin{Def}\label{DefGenPodLoc}
We define the unital $*$-algebra $\mcO_q^{\loc}(S_t^2)$ to be the universal unital $*$-algebra on generators $X_t,Y_t,Z_t^{\pm 1}$ satisfying the same relations as in \eqref{EqCommRelPod1},\eqref{EqCommRelPod2} and \eqref{DefStarQHS}, and with $Z_t^{-1}$ inverse to $Z_t$. 
\end{Def}

This gives an embedding
\[
\mcO_q(S_t^2) \subseteq \mcO_q^{\loc}(S_t^2)
\]
(technically known as an Ore localisation of $\mcO_q(S_t^2)$ at $Z_t$, see e.g.\  \cite{GW04}{Chapter 10}). 

\begin{Prop}
There exists a surjective unital $*$-homomorphism 
\begin{equation}\label{EqFormmsP}
\msP: U_q^{+,-}(\mfsu(2)) \rightarrow \mcO_q^{\loc}(S_t^2)
\end{equation}
such that
\[
(q^{-1}-q)F \mapsto q^{-1/2}Y_t,\qquad (q^{-1}-q)E \mapsto q^{-1/2} Z^{-1}_tX_t,\qquad K \mapsto Z^{-1}_t. 
\]
Moreover, the left $U_q(\mfsu(2))$-module $*$-algebra structure on $\mcO_q(S_t^2)$ given by the infinitesimal right translations 
\[
y \rhd b = (\id\otimes \tau(y,-))\Delta(b),\qquad y \in U_q(\mfsu(2)),b\in \mcO_q(S_t^2)
\] 
extends to $\mcO_q^{\loc}(S_t^2)$, and $\msP$ is then a $*$-homomorphism of $U_q(\mfsu(2))$-module $*$-algebras.
\end{Prop} 

\begin{proof}
We compute directly that $\msP$ is a well-defined unital $*$-homomorphism, which is clearly surjective. Then $\mcO_q^{\loc}(S_t^2)$ is automatically a $U_q(\mfsu(2))$-module $*$-algebra, as the $U_q(\mfsu(2))$-module structure on $U_q^{-,+}(\mfsu(2))$ is inner. Namely, we put
\[
y \blacktriangleright b  = \msP(y_{(+,-;1)})b\msP(S_{-,+}(y_{(-,+;2)})),
\]
and the module $*$-algebra properties are immediate. 

It hence only remains to see that this module structure restricts to $\mcO_q(S_t^2)$ and coincides there with the usual one of $U_q(\mfsu(2))$ by translation. As both structures are module $*$-algebras, we only need to verify that the module action of $E$ and $K$ coincides on the submodule spanned by the generators $1,X_t,Y_t,Z_t$.

But using the matrix presentation $E_t$ in \eqref{EqDefEt}, we see that
\[
y\rhd E_t  = \pi_{1/2}(S(y_{(1)}))E_t \pi_{1/2}(y_{(2)}). 
\]
It is then immediate that $K\rhd  = K \blacktriangleright $ on $E_t$. To see that this also holds for the module action of $E$, we need to check that 
\[
\pi_{1/2}(S(E_{(1)}))E_t \pi_{1/2}(E_{(2)}) =  \begin{pmatrix} 0 & -q^{-1/2} \\ 0 & 0 \end{pmatrix} E_t + \begin{pmatrix} q^{-1} & 0 \\ 0 & q \end{pmatrix} E_t \begin{pmatrix} 0 & q^{1/2} \\ 0 & 0 \end{pmatrix}
\]
equals 
\[
\msP(E_{(+,-;1)})E_t\msP(S_{-,+}(E_{(-,+;2)}))  = (q^{-1}-q)^{-1}q^{-1/2}(Z_t^{-1}X_tE_t - Z_t^{-1}E_tX_t),
\]
which follows from a small computation. 
\end{proof}

\subsection{The regular representation of $\mcU_q(\mfsl(2,\R)_t)$}\label{RegularRep}

In this section, we compute the decomposition of the regular representation 
\[
\pi_{\reg}: \mcU_q(SL(2,\R)_t) \rightarrow \mcB(L^2_q(SL(2,\R)_t)),
\]
where we write 
\[
L^2_q(SL(2,\R)_t) = L^2_q(S_t^2) \otimes L^2_q(SO(2)_t).
\]
We will use for this the specific realisation
\begin{equation}\label{EqInduction2}
L^2_q(SL(2,\R)_t)  = {}_{L^2_q(S_{t}^2)}\Ind(L^2(SO_{q}(2)_t)).
\end{equation}

Indeed, by decomposing $L^2_q(SO(2)_t) \cong \underset{c\in a+\Z}{\bigoplus} p_{[c]}L^2_q(SO(2)_t)$ into its one-dimensional components, we obtain then a corresponding decomposition 
\[
\pi_{\reg} \cong \underset{c\in a+\Z}{\bigoplus}\pi_{\qbra{c}},
\]
where $\pi_{\qbra{c}}$ is the $*$-representation of $\mcU_q(SL(2,\R)_t)$ on $L^2_q(S_t^2)$ given by 
\[
\pi_{\qbra{c}}(a)b = ab,\qquad \pi_{\qbra{c}}(\omega)b = (\omega(b_{(2)}-))_{[c]} b_{(1)},\qquad a,b\in \mcO_q(S_t^2), \omega \in \mcU_q(\mfso(2)_t).
\]
On $L^2_{q,0}(S_t^2)$, we can switch to the corresponding unitary representation of $U_q(\mfsl(2,\R)_t)$. Using the skew-primitivity \eqref{EqSkewPrim} of $iB_t$, we then find the more transparent formula 
\begin{equation}\label{EqConcreteFormulaC}
\pi_{\qbra{c}}(iB_t) b = \tau(b_{(2)},iB_{\qbra{c}}) b_{(1)}.
\end{equation}

The key realisation is now that all these $*$-representations $\pi_{[[c]]}$ can be factorized through embeddings into the global Heisenberg $*$-algebra $\mcH(\mcO_q(S_t^2),U_q(\mfsu(2)))$, acting on $\mcO_q(S_t^2)$ via the unitary $*$-representation
\[
\pi_{\Heis}(a)b = ab,\qquad \pi_{\Heis}(x)b = \tau(b_{(2)},x) b_{(1)},\qquad a,b\in \mcO_q(S_t^2), x\in U_q(\mfsu(2)).
\] 
Indeed, put 
\[
t = \qbra{a},\qquad s = \qbra{c},
\]
and consider the `mixed' infinitesimal Heisenberg $*$-algebra 
\[
\mcH(\mcO_q(S_t^2),U_q(\mfso(2)_s)) \subseteq \mcH(\mcO_q(S_t^2),U_q(\mfsu(2)))
\] 
as in \eqref{EqMixedHeis}.  Upon realizing that we can write
\[
\Delta(iB_s) = iB_t \otimes 1 + K\otimes (iB_s -[a]),
\] 
it follows immediately from the interchange relations \eqref{EqCommRelMixedHeis} that we obtain an isomorphism of $*$-algebras 
\begin{equation}\label{EqMapIntoHeis}
\theta_{s,t}: U_q(\mfsl(2,\R)_t)  \rightarrow \mcH(\mcO_q(S_t^2),U_q(\mfso(2)_s)),\qquad b \mapsto b,\quad iB_t \mapsto iB_s,\qquad b\in \mcO_q(S_t^2), 
\end{equation}
and \eqref{EqConcreteFormulaC} tells us that 
\begin{equation}\label{EqCompOfReps}
\pi_{\qbra{c}} = \pi_{\Heis}\circ \theta_{s,t}. 
\end{equation}
Below, we will also use the `mixed Casimir'
\begin{equation}\label{EqMixedCas}
\Omega_{s,t} := \theta_{s,t}(\Omega_t) =   iq^{-1}X_t +(q-q^{-1})iZ_tB_s-iqY_t\in \mcH(\mcO_q(S_t^2),U_q(\mfso(2))_t)
\end{equation}
that one obtains under the above isomorphism, with $\Omega_t$ the ordinary Casimir element as defined in \eqref{EqCasimir}.

We now show that $\pi_{\Heis}$ has a compatible realisation using instead the pre-Hilbert spaces $V_{\pm},W_{\pm}$ introduced in \eqref{EqPreHilb1} and \eqref{EqPreHilb2}. 

\begin{Rem}
We note that a detailed study of $*$-representations of $\mcH(\mcO_q(S_t^2),U_q(\mfsu(2)))$ was made in \cite{SW07}. Although we will not use any results directly, the use of the localisation of $\mcO_q(S_t^2)$ is exactly what allows to have a `decoupling' of the Heisenberg algebra (= crossed product algebra) $\mcH(\mcO_q^{\loc}(S_t^2),U_q(\mfsu(2)))$ as in \cite{SW07}*{Section 5} (see also \cite{Fio02}). This will be what is behind the scenes in the following computations. 
\end{Rem} 
The following theorem appears (in a somewhat different formulation) in \cite{SW07}*{Section 6}. We present the proof (in a slightly more conceptual way) for the benefit of the reader. Recall the notation \eqref{EqPreHilb2}.

\begin{Theorem}\label{TheoEqTranspo}
There exist $*$-representations $\widetilde{\pi}_{\pm}$ of $\mcH(\mcO_q(S_t^2),U_q(\mfsu(2)))$ on $W_{\pm}$ such that the following conditions hold: 
\begin{itemize}
\item $\widetilde{\pi}_{\pm}$ coincide with the usual GNS-action on $\mcO_q(S_t^2)$, so 
\begin{equation}\label{EqGNSRep}
\widetilde{\pi}_{\pm}(b) = \pi_{\pm}(b)\otimes 1,\qquad b\in \mcO_q(S_t^2),
\end{equation}
with  $\pi_{\pm}$ defined in Proposition \ref{PropRepPod},
\item on the generators $E,F,K\in U_q(\mfsu(2))$ we have 
\[
\widetilde{\pi}_{\pm}(K) = Z_t^{-1}\otimes Z_t, 
\]
\[
(q^{-1}-q)\widetilde{\pi}_{\pm}(E) = q^{-1/2} Z_t^{-1}X_t\otimes 1 - q^{1/2} Z_t^{-1}\otimes Y_t
\]
\[
(q^{-1}-q)\widetilde{\pi}_{\pm}(F) = q^{-1/2} Y_t\otimes Z_t^{-1} - q^{1/2}1\otimes X_tZ_t^{-1} 
\]
and
\item for all $\xi\in L^2_{q,0}(S_t^2),\eta\in W^{\pm}$ and $x\in \mcH(\mcO_q(S_t^2),U_q(\mfsu(2)))$ one has that
\begin{equation}\label{EqCorrHeisStar}
\langle \widetilde{\pi}_{\pm}(x)\eta,\xi\rangle = \langle \eta,\pi_{\Heis}(x^*)\xi\rangle. 
\end{equation}
\end{itemize}
\end{Theorem} 
\begin{proof}
We have a linear $*$-anti-isomorphism
\[
R_{-,+}: U_q^{-,+}(\mfsu(2))\rightarrow U_q^{+,-}(\mfsu(2)),\qquad K \mapsto K^{-1},\quad E \mapsto -qK^{-1}E,\quad F \mapsto -q^{-1}FK. 
\]
Denoting by $\overline{R}_{-,+}$ the anti-linear anti-isomorphism defined by the same formula on generators, we can then check that the action of the generators in the second point of the theorem is given by  
\[
\widetilde{\pi}_{\pm}(x) = \msP(x_{(+,-;1)})\otimes \msP(\overline{R}_{-,+}(x_{(-,+;2)})^*),
\]
where $\msP$ is the map introduced in \eqref{EqFormmsP}.
 So, we at least get a $*$-representation of $U_q(\mfsu(2))$. 

To see that it extends to a $*$-representation of the Heisenberg $*$-algebra via \eqref{EqGNSRep}, it is sufficient to check that \eqref{EqCorrHeisStar} is satisfied for $x\in U_q(\mfsu(2))$. Recall that $T = \frac{q^{-1}-q}{\qdif{a}}|Z_t| \in \mcB(\mcH_+ \oplus \mcH_-)$. Let us write again explicitly $\xi_T \in L^2_q(S_t^2)$ for the GNS-vector of $\mcO_q(S_t^2)$ with respect to $\Phi= \Tr(-T)$. We then compute that for $b\in \mcO_q(S_t^2)$, using the notation \eqref{EqDoubleInd}, 
\begin{eqnarray*} 
\langle \widetilde{\pi}_{\pm}(x)\xi_{m,n}^{\pm},(b\otimes 1)\xi_T\rangle &=&\langle (b^*\otimes 1)\widetilde{\pi}_{\pm}(x)\xi_{m,n}^{\pm},\xi_T\rangle \\
&=& \langle (b^*\msP(x_{(+,-;1)})\otimes \msP(\overline{R}_{-,+}(x_{(-,+;2)})^*))\xi_{m,n}^{\pm},\xi_T\rangle \\
&=& \langle  (b^*\msP(x_{(+,-;1)})\xi_m^{\pm} ,T^{1/2} \msP(R_{-,+}(x_{(-,+;2)})^*))\xi_n^{\pm}\rangle.
\end{eqnarray*}
Now it is easily verified that $T^{1/2}\msP(R_{-,+}(x)) = \msP(S_{-,+}(x))T^{1/2}$. Hence 
\begin{eqnarray*}
\langle \widetilde{\pi}_{\pm}(x)\xi_{m,n}^{\pm},(b\otimes 1)\xi_T\rangle &=& \langle  (b^*\msP(x_{(+,-;1)})\xi_m^{\pm} , \msP(S_{-,+}(x_{(-,+;2)}^*)))T^{1/2}\xi_n^{\pm}\rangle \\
 &=& \langle  \xi_m^{\pm} ,\msP(x_{(+,-;1)}^*) b \msP(S_{-,+}(x_{(-,+;2)}^*)))T^{1/2}\xi_n^{\pm}\rangle\\
&=& \langle  \xi_m^{\pm} ,(x^* \rhd b)T^{1/2}\xi_n^{\pm}\rangle\\ 
&=& \langle \xi_{m,n}^{\pm},((x^*\rhd b)\otimes 1)\xi_T\rangle\\
&=& \langle \xi_{m,n}^{\pm},\pi_{\Heis}(x^*)(b\otimes 1)\xi_T\rangle.
\end{eqnarray*} 
Since $L^2_{q,0}(S_t^2) \subseteq L^2_q(S_t^2)$ consists precisely of elements of the form $(b\otimes 1)\xi_T$, we are done. 
\end{proof}

Using now \eqref{EqMapIntoHeis}, we obtain a $*$-representation 
\[
\widetilde{\pi}_s^{\pm} := \widetilde{\pi}_{\pm}\circ \theta_{s,t}
\]
of $U_q(\mfsl(2,\R)_t)$ on $W^{\pm}$, given more directly as  
\begin{equation}\label{EqDefOfRepTildePi}
\widetilde{\pi}_s^{\pm}(b) = \pi_{\pm}(b)\otimes 1,\qquad \widetilde{\pi}_s^{\pm}(iB_t) = \widetilde{\pi}_{\pm} (iB_s),\qquad b\in \mcO_q(S_t^2). 
\end{equation}
By Theorem \ref{TheoEqTranspo}, these representations satisfy
\begin{equation}\label{EqUnderTheAdjointCorr}
\langle \widetilde{\pi}_s^{\pm}(x)\xi,\eta\rangle = \langle \xi,\pi_s(x^*)\eta\rangle,\qquad \xi \in W^{\pm},\eta \in L^2_{q,0}(S_t^2),x\in U_q(\mfsl(2,\R)_t).
\end{equation}
Concretely, we have 
\begin{equation}\label{EqFormulaWPZ}
\widetilde{\pi}_s^{\pm}(Z_t)\xi_{p,k}^{\pm} =\pm q^{2p\mp a+1} \xi_{p,k}^{\pm},
\end{equation}
\begin{equation}\label{EqFormulaWPY}
\widetilde{\pi}_s^{\pm}(Y_t)\xi_{p,k}^{\pm}=  \sqrt{(1-q^{2p+2})(q^{\mp 2a+2p+2}+1)}\xi_{p+1,k}^{\pm},
\end{equation}
\begin{equation}\label{EqFormulaWPX}
\widetilde{\pi}_s^{\pm}(X_t) \xi_{p,k}^{\pm} = \sqrt{(1-q^{2p})(q^{\mp 2a+2p}+1)}\xi_{p-1,k}^{\pm},
\end{equation}
\begin{equation}\label{EqFormulaWPB}
(q-q^{-1})i\widetilde{\pi}_s^{\pm}(B_t) \xi_{p,k}^{\pm} = (i(q^{-1}Y_t-qX_t)Z_t^{-1} \otimes 1 +Z_t^{-1}\otimes (i(Y_t-X_t)+sZ_t))\xi_{p,k}^{\pm},
\end{equation}
where on the right hand side of \eqref{EqFormulaWPB} the elements $Y_t,X_t,Z_t^{\pm1}$ act via $\pi_{\pm}$ on both legs. 

Consider now the mixed Casimir of \eqref{EqMixedCas}. Then we see from \eqref{EqDefOfRepTildePi}, \eqref{EqFormulaWPB} and the commutation relations in $\mcO_q^{\loc}(S_t^2)$ that 
\begin{equation}\label{EqOmForm}
\widetilde{\pi}_s^{\pm}(\Omega_t) = \widetilde{\pi}_{\pm}(\Omega_{s,t}) =  1\otimes \pi_{\pm}(i(Y_t-X_t)+ sZ_t). 
\end{equation}
In the following, we write 
\begin{equation}\label{EqOmTilde}
\widetilde{\Omega}_{s,t} := i(Y_t-X_t)+ sZ_t.
\end{equation}
Note that this is a `balanced' (and self-adjoint) version of the operators $A_{c}$ as defined in \cite{DCDT24a}. We have
\begin{multline*}
\pi_{\pm}(\widetilde{\Omega}_{\qbra{c},\qbra{a}})\xi_p^{\pm} =  i\sqrt{(1-q^{2p+2})(1+q^{\mp 2a+2p+2})}\xi_{p+1}^{\pm} \pm  \qbra{c} q^{2p\mp a+1} \xi_{p}^{\pm} -i  \sqrt{(1-q^{2p})(1+q^{\mp 2a+2p})}\xi_{p-1}^{\pm},
\end{multline*}
which defines a bounded, self-adjoint operator.

We will need some information on the spectrum of $\pi_{\pm}(\widetilde{\Omega}_{\qbra{c},\qbra{a}})$. For this, we first recall that the \emph{Al-Salam Chihara polynomials} $Q_p(\lambda) = Q_p(\lambda;\alpha,\beta\mid q^2)$ (depending on $\alpha,\beta\in \C$ and $q>0$) are the unique polynomials in $\lambda$ satisfying the recurrence 
\begin{equation}\label{EqAlSalam}
Q_0(\lambda) = 1,\qquad 2\lambda Q_p(\lambda) = Q_{p+1}(\lambda) +(\alpha + \beta)q^{2p}Q_p(\lambda) +(1-q^{2p})(1-\alpha\beta q^{2p-2})Q_{p-1}(\lambda).
\end{equation}

\begin{Lem}
Let
\[
\xi^{\pm}(\lambda) = \sum_{p=0}^{\infty} A_p^{\pm}(\lambda)\xi_p^{\pm} \in \prod_{n\geq 0} \C\xi_{n}^{\pm}
\]
be a (suitably normalized) formal eigenvector for $\pi_{\pm}(\widetilde{\Omega}_{\qbra{c},\qbra{a}})$. Then 
\begin{equation}\label{EqFormEigVect}
A_p^{\pm}(\lambda) = \frac{Q_p(-\lambda/2;\mp q^{\mp a +c+1},\pm q^{\mp a -c +1}\mid q^2)}{i^p(q^2;q^2)_p^{1/2}(-q^{\mp 2a+2};q^2)_p^{1/2}}.
\end{equation}
\end{Lem}
Here we use the standard notation for $q^2$-Pochhammer symbols, so 
\[
(\alpha;q^2)_p = (1-\alpha)(1-q^2\alpha)\ldots (1-q^{2p-2}\alpha),\qquad \alpha\in \C,p\in \Z_{\geq0}. 
\]
\begin{proof}
Using the explicit expression for $\pi_{\pm}(\widetilde{\Omega}_{\qbra{c},\qbra{a}})$, we obtain the relation 
\begin{multline}
\lambda A_p^{\pm}(\lambda) = i\sqrt{(1-q^{2p})(1+q^{\mp 2a+2p})} A_{p-1}^{\pm}(\lambda) \pm \qbra{c} q^{2p\mp a+ 1}A_p^{\pm}(\lambda) \\ - i \sqrt{(1-q^{2p+2})(1+q^{\mp 2a+2p+2})}A_{p+1}^{\pm}(\lambda).
\end{multline}
Put then 
\[
B_p^{\pm}(\lambda) =i^p(q^2;q^2)_p^{1/2}(-q^{\mp 2a+2};q^2)_p^{1/2} A_p^{\pm}(\lambda).
\]
Using the (defining) identities 
\[
(q^2;q^2)_p = (q^2;q^2)_{p-1}(1-q^{2p}),\quad (-q^{\mp 2a + 2};q^2)_p = (-q^{\pm2a + 2};q^2)_{p-1}(1+q^{\pm 2a +2p})
\]
for $q^2$-Pochhammer symbols, we see immediately that $B_p(\lambda)$ satisfies the recurrence relation 
\[
-\lambda B_p^{\pm}(\lambda) = B_{p+1}^{\pm}(\lambda)  \mp \qbra{c}q^{2p\mp a +1}B_p^{\pm}(\lambda) + (1-q^{2p})(1+q^{\mp 2a +2p})B_{p-1}^{\pm}(\lambda). 
\]
Comparing with \eqref{EqAlSalam}, we can conclude \eqref{EqFormEigVect}.
\end{proof}

We note in passing that the Al-Salam-Chihara polynomials allow for a closed form expression: using standard notation for basic hypergeometric functions \cite{GR04}, one can write (in the case of our specific parameters)
\[
Q_p(\cos(\theta);\mp q^{\mp a +c+1},\pm q^{\mp a -c +1}\mid q^2) := \frac{(-q^{\mp 2a+2};q^2)_p}{(\mp 1)^p q^{p(c\mp a + 1)}}\rphis{3}{2}{q^{-2p},\mp q^{\mp a +c+1}e^{i\theta},\mp q^{\mp a +c+1}e^{-i\theta}}{-q^{\mp 2a+2},0}{q^2,q^2}, 
\]
see again \cite{KLS10}*{Section 14.8}.

Now from \cite[Section 3.3 and Section 3.4]{AI84}, we obtain the following. 
\begin{Theorem}\label{TheoTwistedCasimir}
The spectrum of $\pi_{\pm}(\widetilde{\Omega}_{\qbra{c},\qbra{a}})$ is multiplicity-free, with the continuous spectrum located on $[-2,2]$ and the point-spectrum determined as the finite set $J_{\qbra{c},\qbra{a}}^{\pm} \subseteq (-\infty,-2) \cup (2,\infty)$ with 
\[
J_{\qbra{c},\qbra{a}}^{\pm} = J_{\qbra{c},\qbra{a}}^{\pm}(+) \cup J_{\qbra{c},\qbra{a}}^{\pm}(-)
\]
with\footnote{In \cite[Section 3.3]{AI84}, the respective cases $k =  \frac{\pm(a+c)-1}{2}$ and $k =  \frac{\pm(a-c)-1}{2}$ (when these are positive integers) need to be excluded, as the points $\pm 2$ do not belong to the point spectrum, as shown there. This is also clear from the value of the discrete weight below, which vanishes in these points.} 
\begin{equation}\label{EqJValPos}
J_{\qbra{c},\qbra{a}}^{\pm}(+) = \{\qdif{\pm(a-c)-2k-1}\mid k\in  \Z_{\geq0} \cap \Z_{< \frac{\pm(a-c)-1}{2}}\},
\end{equation}
\begin{equation}\label{EqJValNeg}
J_{\qbra{c},\qbra{a}}^{\pm}(-) = \{-\qdif{\pm(a+c)-2k-1}\mid k\in \Z_{\geq0} \cap \Z_{<  \frac{\pm(a+c)-1}{2}}\}.
\end{equation}

Moreover, the spectral measure $\mu_{\qbra{c},\qbra{a}}^{\pm}$ with respect to $\xi_0^{\pm}$ is of the form 
\[
\rd\mu_{\qbra{c},\qbra{a}}^{\pm}(\lambda) = g_{\qbra{c},\qbra{a}}^{\pm}(\lambda) \rd \lambda + \sum_{\lambda' \in J_{\qbra{c},\qbra{a}}^{\pm}} w_{\qbra{c},\qbra{a}}^{\pm}(\lambda') \delta_{\lambda'}(\lambda),
\]
where
\[
g_{\qbra{c},\qbra{a}}^{\pm}(e^{i\theta} + e^{-i\theta}) = \frac{1}{4\pi |\sin(\theta)|} \frac{(e^{2i\theta},e^{-2i\theta},-q^{\mp 2a + 2},q^2;q^2)_{\infty}}{(q^{\mp(a-c)+1}e^{i\theta},-q^{\mp(a+c)+1}e^{i\theta},q^{\mp(a-c)+1}e^{-i\theta},-q^{\mp(a+c)+1}e^{-i\theta};q^2)_{\infty}}
\]
and with 
\[
w_{\qbra{c},\qbra{a}}^{\pm}(\qdif{\pm(a-c)-2k-1}) = q^{2k}(1-q^{\pm 2(a-c)-4k-2})\frac{(-q^{\pm2a -2k};q^2)_k(q^{\pm 2(a-c)-2k};q^2)_{\infty}}{(q^2;q^2)_k (-q^{\mp 2c};q^2)_{\infty}(-q^{2\pm 2c};q^2)_k},
\]
\[
w_{\qbra{c},\qbra{a}}^{\pm}(-\qdif{\pm(a+c)-2k-1}) = q^{2k}(1- q^{\pm 2(a+c) -4k-2}) \frac{(-q^{\pm2a -2k};q^2)_k(q^{\pm 2(a+c)-2k};q^2)_{\infty}}{(q^2;q^2)_k (-q^{\pm 2c};q^2)_{\infty}(-q^{2 \mp 2c};q^2)_{k}}.
\]
\end{Theorem}

\begin{Rem}\label{RemDisj}
Implicit in the above theorem is the fact that $J_{\qbra{c},\qbra{a}}^+$ and $J_{\qbra{c},\qbra{a}}^-$ are disjoint, which follows from an elementary calculation using that $\qdif{x} = \qdif{y}$ if and only if $x= y$ or $x=-y$. 
\end{Rem} 

We can now prove: 
\begin{Theorem}\label{TheoDecompPrinc}
Assume that $t = \qbra{a}$, with $4a$ not an odd integer. Then the regular representation $\pi_{\reg}$ of quantum $SL(2,\R)$ decomposes as 
\begin{equation}\label{EqDecompFund}
L^2_q(SL(2,\R)_t) \cong \bigoplus_{\N}\left(\left( \oplus_{\pm} \int_{[-2,2]} \mcL_{\lambda}^{\pm} \rd \lambda\right)\oplus \bigoplus_{n\geq 2} \mcD_n^+ \oplus \bigoplus_{n\geq 2} \mcD_n^- \oplus \bigoplus_{n\geq 2} \mcE_n^+ \oplus \bigoplus_{n\geq 2} \mcE_n^-\right),
\end{equation}
where the direct integral decomposition is with respect to Lebesgue measure.

In case $4a$ is an odd integer, the indexing for the exceptional discrete series $\mcE_n^+$ and $\mcE_n^-$ needs to start at $1$. 
\end{Theorem}
We do not have a good conceptual reason why the $\mcE_1^{\pm}$ dip into the decomposition when $4a \in 2\Z+1$. We only remark that these are exactly the values where $\mcE_{1}^+$ and $\mcE_{2}^+$ (and $\mcE_{1}^-$ and $\mcE_{2}^-$)  share the same value of the Casimir.\footnote{Also, should in this particular case $\mcE_1^{\pm 1}$ be called mock mock exceptional discrete series representations?}
\begin{proof}
As a unitary $\mcO_q(SU(2))$-comodule, we have 
\[
L^2_{q,0}(S_t^2) \cong \oplus_{n\in \Z_{\geq0}} V_{n},
\]
where we recall that $V_n$ denotes the spin $n$-representation of $U_q(\mfsu(2))$. Hence by \cite{Koo93}*{Theorem 4.3}, the action of $iB_{\qbra{a+n}}$ on $L^2_{q,0}(S_t^2)$ has spectrum $\{[a+n+2k]\mid k\in \Z\}$, and the eigenvalues of $\pi_{\qbra{a+n}}(iB_t)$ are
\begin{itemize}
\item of the form $\{[a+2k]\mid k \in \Z\}$ if $n$ is even, and
\item of the form $\{[a+2k+1]\mid k \in \Z\}$ if $n$ is odd.
\end{itemize}
Now from the decomposition of the Casimir in Theorem \ref{TheoTwistedCasimir} and the list of values for the Casimir and spectrum of $iB_t$ for the irreducible representations  in Section \eqref{SecIrrRepSL2}, we  see that $\pi_{\qbra{a+n}}$ must contain $\int_{[-2,2]} \mcL_{\lambda}^+ \rd \lambda$ in case $n$ is even, and $\int_{[-2,2]} \mcL_{\lambda}^- \rd \lambda$ in case $n$ is odd. This accounts for the continuous part of the decomposition \eqref{EqDecompFund}.

For the discrete components arising in \eqref{EqDecompFund}, we argue as follows. 

First note that the identities \eqref{EqOmForm}, \eqref{EqOmTilde} and \eqref{EqUnderTheAdjointCorr}, together with Theorem \ref{TheoTwistedCasimir}, imply that the discrete spectrum (= set of eigenvalues) of $\pi_{\qbra{a+n}}(\Omega_t)$ equals the set $J_{\qbra{a+n},\qbra{a}}^+ \cup J_{\qbra{a+n},\qbra{a}}^-$. As $J_{\qbra{a+n},\qbra{a}}^+$ and $J_{\qbra{a+n},\qbra{a}}^-$ are disjoint (Remark \ref{RemDisj}), we can choose for each such eigenvalue $\lambda$  a unique sign $\pm = \pm(\lambda)$ such that $\lambda \in J_{\qbra{a+n},\qbra{a}}^{\pm}$. In the following, we fix $\lambda \in J_{\qbra{a+n},\qbra{a}}^{\pm}$, and we let $q_{\lambda}$, resp.\ $p_{\lambda}$ be the spectral projection of $\pi_{\qbra{a+n}}(\Omega_t)$, resp.\  $\pi_{\pm}(\widetilde{\Omega}_{\qbra{a+n},\qbra{a}})$ at eigenvalue $\lambda$. Then under the identification \eqref{EqIdentificationPodSph}, we get 
\begin{equation}\label{EqIdentificationValCas}
q_{\lambda}L^2_q(S_t^2) \cong \Hsp_{\pm} \otimes p_{\lambda}\Hsp_{\pm},
\end{equation}
where $p_{\lambda}\Hsp_{\pm}$ is one-dimensional.

Now consider the restriction $\pi_{\qbra{a+n}}^{(\lambda)}$ of $\pi_{\qbra{a+n}}$ to $q_{\lambda}L^2_q(S_t^2)$. It will be a $*$-representation of $\mcU_q(SL(2,\R)_t)$ such that 
\begin{itemize}
\item on the corresponding algebraic part, the Casimir $\Omega_t$ takes on the constant value $\lambda$, 
\item on the corresponding algebraic part, the element $iB_t$ has eigenvalues lying in the set $\{[a+n+2k]\mid k \in \Z\}$, and
\item the element $Z_t$ is either strictly positive or strictly negative, and more particularly of constant sign $\pm(\lambda)$, as one sees by using \eqref{EqIdentificationValCas}, \eqref{EqUnderTheAdjointCorr} and the concrete formula \eqref{EqFormulaWPZ}.
\end{itemize}
Clearly, these properties are inherited by almost all irreducible representations appearing in a direct integral decomposition of $\pi_{\qbra{a+n}}^{(\lambda)}$ (cf.\ Remark \ref{RemType1}).

Now by Proposition \ref{PropRepPod} and Theorem \ref{TheoBranching}, the last property above guarantees that only discrete series and exceptional discrete series representations can appear. Moreover, the sign of $Z_t$ then dictates which ones can precisely appear: 
\begin{itemize}
\item Only $\mcD_{n'}^-$ and $\mcE_{n'}^+$ can appear if $\pm(\lambda) = +$.
\item Only $\mcD_{n'}^+$ and $\mcE_{n'}^-$ can appear if $\pm(\lambda) = -$.
\end{itemize}
Since moreover the index $n'$ can be recovered from the knowledge of the value of the Casimir operator and the parity of the operator $iB_t$ (again by direct inspection of the list in Section \ref{SecIrrRepSL2}), we conclude that $\pi_{\qbra{a+n}}^{(\lambda)}$ must be, up to multiplicity, a single discrete series (or exceptional discrete series) representation (with the mock ones at the moment still under consideration).

To conclude the theorem, we must now determine which (exceptional) discrete series representations can arise when $n$ varies. We will focus on the case of $\mcE_{n'}^{\pm}$ for $n' \in \Z_{>0}$ (the case $\mcD_{n'}^{\pm}$ is similar, but easier). 

If $n'\in \Z_{>0}$, we have by the above that $\mcE_{n'}^+$ appears in $\pi_{\qbra{a+n}}$ if and only if it appears in $\pi_{\qbra{a+n}}^{(\lambda)}$ for some $\lambda \in J_{\qbra{a+n},\qbra{a}}^+(-)$. Considering the precise form for these elements given by \eqref{EqJValNeg}, we see that this will occur if and only if there exists $k\in \Z_{\geq0}$ with 
\begin{equation}\label{EqEqualAbsValue}
|2a + n-2k -1| = |2a +\lfloor -2a\rfloor + n'-1|\qquad \textrm{ and }\qquad  2k < 2a + n-1. 
\end{equation}
Now rewriting the latter equation as $-2a+1 < n-2k$, we see that either there will be infinitely $n$ allowing such a $k$, or none at all. Moreover, the left hand side of \eqref{EqEqualAbsValue} gives that either 
\[
n-2k = \lfloor -2a\rfloor + n'\qquad \textrm{or}\qquad n-2k = -4a - \lfloor -2a \rfloor - n' + 2.
\]
So for there to exist one (and hence infinitely many) couples $(n,k)$ offering a solution, at least one of the following properties must hold: 
\begin{enumerate}
\item We have $-2a + 1 < \lfloor -2a\rfloor + n'$. Clearly this holds if and only if $n' \geq 2$. 
\item We have that $4a$ is an integer, and $-2a + 1< -4a - \lfloor -2a\rfloor -n' +2$, which can be rewritten as $n' < 1-2a - \lfloor -2a \rfloor$. Clearly this holds if and only if $4a$ is an odd integer and $n' = 1$.  
\end{enumerate} 
This gives us the correct count of $\mcE_{n'}^+$'s in Theorem \ref{TheoDecompPrinc}. 

The case of $\mcE_{n'}^-$'s is similar: if $n'\in \Z_{>0}$, we have that $\mcE_{n'}^-$ appears in $\pi_{\qbra{a+n}}$ if and only if it appears in $\pi_{\qbra{a+n}}^{(\lambda)}$ for some $\lambda \in J_{\qbra{a+n},\qbra{a}}^-(-)$, so this will occur if and only if there exists $k\in \Z_{\geq0}$ with 
\begin{equation}\label{EqEqualAbsValue2}
|2a + n+2k +1| = |2a -\lfloor 2a\rfloor - n'+1|\qquad \textrm{ and }\qquad  2k < -2a - n-1 \iff n+2k < -2a-1.
\end{equation}
Again, we see that there then will either be infinitely $n$ allowing such a $k$, or none at all. Decoupling again the left hand side of  \eqref{EqEqualAbsValue2} gives that either 
\[
n+2k =- \lfloor 2a\rfloor - n'\qquad \textrm{or}\qquad n+2k = -4a + \lfloor 2a \rfloor + n' -2.
\]
So for there to exist one (and hence infinitely many) couples $(n,k)$ offering a solution, at least one of the following properties must hold: 
\begin{enumerate}
\item We have $ -\lfloor 2a\rfloor - n' < -2a-1$. Clearly this holds if and only if $n' \geq 2$. 
\item We have that $4a$ is an integer, and $-4a + \lfloor 2a\rfloor +n' -2<-2a-1$, which can be rewritten as $n' < 2a -\lfloor 2a\rfloor +1$. Clearly this holds if and only if $4a$ is an odd integer and $n' = 1$.  
\end{enumerate} 
This gives us the correct count of $\mcE_{n'}^-$'s in Theorem \ref{TheoDecompPrinc}. 
\end{proof}

\begin{Rem}
With a little more effort, one can show that $\pi_{\qbra{a+n}}$ is multiplicity-free (i.e.\ its space of intertwiners is commutative). Let us sketch the argument.

Note first that we can uniquely extend $\pi_{\qbra{a+n}}(iB_t)$ to a self-adjoint operator (with discrete spectrum) on $L^2_q(S_t^2)$, which we will write $i\mathbf{B}_{\qbra{a+n}}$. Since $L^2_{q,0}(S_t^2)$ is a core for this operator, it follows from Theorem \ref{TheoEqTranspo} that, under the identification 
\begin{equation}\label{EqIdentificationPodSph}
L^2_q(S_t^2) \cong \oplus_{\pm} (\Hsp_{\pm} \otimes \Hsp_{\pm}),
\end{equation}
the transported operator $i\mathbf{B}_{\qbra{a+n}}$ will have $W_{\pm}$ in its domain, where it acts as $\widetilde{\pi}_{\pm}(iB_{\qbra{a+n}})$. 

Now the (bounded) selfadjoint operator $\pi_{\qbra{a+n}}(\Omega_t)$ strongly commutes with $i\mathbf{B}_{\qbra{a+n}}$. So for $\lambda$ in the discrete spectrum of $\pi_{\qbra{a+n}}(\Omega_t)$, we see from the specific form \eqref{EqOmForm} that under the identification \eqref{EqIdentificationValCas}, we get
\begin{equation}\label{EqProjB}
(q-q^{-1})(i\mathbf{B}_{\qbra{a+n}})(\xi\otimes p_{\lambda}\eta) = (iq^{-1}Y_tZ_t^{-1} -iqX_tZ_t^{-1} + \lambda Z_t^{-1})\xi\otimes p_{\lambda}\eta,\qquad \xi,\eta \in V^{\pm}. 
\end{equation}
As this is a self-adjoint operator restricting to a tridiagonal operator on an orthonormal basis within its domain, it follows that its spectrum must be simple. Since all eigenvalues of  $i\mathbf{B}_{\qbra{a+n}}$ differ by an even integer, it follows that indeed $\pi_{\qbra{a+n}}^{(\lambda)}$ must be a simple representation. This takes care of the claim for the discrete spectrum. 

For the continuous spectrum one can follow a similar approach, but one needs to be a bit more careful: now any value  in $(-2,2)$ will appear in the spectrum of $\pi_{\pm}(\widetilde{\Omega}_{\qbra{a+n}},\qbra{a})$ for both choices of sign, and we can no longer decompose $i\mathbf{B}_{\qbra{a+n}}$ as a direct sum of self-adjoint operators over the two Hilbert spaces $\int_{[-2,2]}^{\oplus} \Hsp_{\pm} \otimes p_{\lambda}\Hsp_{\pm} \rd \lambda$. However, this \emph{does} tell us that the multiplicity function takes on only the values $1$ and $2$ on the continuous spectrum. Since however $Z_t$ needs to assume both positive and negative values in a principal series representation by Theorem \ref{TheoBranching}, we can conclude that the multiplicity function needs to be constant $1$ also on the continuous part.
\end{Rem}

\section{Outlook}

In this paper, we have explained in general how, given a unitary Doi-Koppinen datum of coideal type $(A,B,C)$, one can construct canonical unitary Doi-Koppinen representations $L^2(B) \in {}_B\Rep^A$ and $L^2(\mcI) \in  {}_A\Rep^C$,  with $\mcI$ is the restricted dual of $C$. We have then shown how these unitary representations lead to respective induction functors 
\[
{}_{L^2(B)}\Ind: \Rep^C \rightarrow {}_B\Rep^C,\qquad \Ind^{L^2(\mcI)}: {}_B\Rep \rightarrow {}_B\Rep^C.  
\]
We have moreover introduced the regular representation 
\[
L^2(B) \otimes L^2(\mcI) = {}_{L^2(B)}\Ind(L^2(\mcI)) \cong \Ind^{L^2(\mcI)}(L^2(B)), 
\]
of $(A,B,C)$, and constructed as well a canonical representation on $L^2(A)$. We have then computed how these operations behave in the case of the unitary DK-datum 
\[
(\mcO_q(SU(2)),\mcO_q(S_t^2),\mcO_q(SO(2)_t)),
\]
whose unitary DK-representations can be interpreted as the unitary representations of quantum $SL(2,\R)$. In particular, we have computed how the regular representation of quantum $SL(2,\R)$ decomposes into irreducibles. 

Here are some questions that were left unanswered, and that will be left for a future occasion. 

\begin{itemize}
\item One can also consider the `second' canonical representation of quantum $SL(2,\R)$ on $L^2_q(SU(2))$. At the moment, it is unclear if this representation will be unitarily equivalent to the regular representation. As of yet, we have not attempted to compute the decomposition of this representation into irreducibles. This would at least lead directly to an answer to the above question, but would not really be illuminating from a conceptual standpoint. 
\item Given that quantum $SL(2,\R)$ posesses a natural positive (relatively invariant) functional $\varphi_t$ as constructed through \eqref{EqCanFunctDD}, it is natural to consider a Plancherel formula for $\varphi_t$ with respect to the decomposition of the regular representation into irreducibles. Again, further considerations will be necessary to obtain concrete results in this direction. 
\item It would be interesting to understand in more detail the limit behaviour for $q \rightarrow 1$ of (the objects of and constructions for) the representation theory of $\mcU_q(SL(2,\R)_t)$ from the analytic point of view (in the formal setting, this was briefly considered in \cite{DCDT24a}). Also the dependence of analytic features on the parameter $t$ (as e.g.\ in Theorem \ref{TheoDecompPrinc}) should be understood better. 
\end{itemize}

\emph{Acknowledgements}: This research was funded by the FWO grant G032919N. The author thanks J.R. Dzokou Talla for discussions on the results in Section 5, and the referees for their detailed and insightful comments, which have improved the results in this paper.

\end{document}